\documentclass{article}
\usepackage{graphicx}
\usepackage{bigints}

\usepackage[utf8]{inputenc}
\usepackage[margin=1in]{geometry}
\usepackage{color}
\usepackage{enumitem}
\usepackage{authblk}
\usepackage{esint}
\usepackage{lineno}

\usepackage{hyperref}
\usepackage{hypcap}

\usepackage{algorithm}
\usepackage{algpseudocode}

\usepackage{float}
\usepackage{caption}
\usepackage{subcaption}
\usepackage{setspace}
\usepackage{tcolorbox}

\usepackage{amsmath,amsfonts,amssymb,mathtools, amsthm}
\usepackage[english]{babel}

\newtheorem{defi}{Definition}

\newtheorem{remark}{Remark}
\numberwithin{equation}{section}

\usepackage{authblk}

\title{A fast scheme for the homogeneous Boltzmann equation based on lifting and tensor train approximation}

\author{Kun Huang  \thanks{Department of Mathematics, Virginia Tech, Blacksburg, VA 24061 U.S.A. (\textsf{kunhuang@vt.edu})} \, 
Yingda Cheng \thanks{Department of Mathematics, Virginia Tech, Blacksburg, VA 24061 U.S.A. (\textsf{yingda@vt.edu}). Research is supported by DOE grant DE-SC0023164, AFOSR grant FA9550-25-1-0154 and Virginia Tech.} \, Irene M. Gamba \thanks{Department of Mathematics and Oden Institute for Computational Sciences and Engineering, University of Texas at Austin (\textsf{gamba@math.utexas.edu}). Research is supported by NSF grant DMS 2408263.}}

\date{}
\begin{document}
\maketitle

\begin{abstract} 
We propose a fast deterministic scheme for the space-homogeneous Boltzmann equation that exploits the low-rank structure of the velocity distribution. This paper consists of two independent contributions. The first is a \emph{lifting-projection (LP) scheme}, inspired by the approach in the recent theoretical breakthroughs \cite{guillen2025landau, imbert2026monotonicity, guillen2025landau2} on the well-posedness of the Landau and Boltzmann equations. In particular, the approach lifts the nonlinear 3D Boltzmann equation  to the 6D linear Kac master equation, advanced over a single time step, and projected back to its marginal in 3D. The second contribution is a \emph{low-rank tensor method} for evaluating the collision operator, in which the lifted solution is represented in tensor train (TT) format and computed via a TT cross approximation algorithm with interpolation, complemented by a TT-friendly conservation correction that enforces conservation of mass, momentum, and energy. When the solution is low-rank in velocity, the method scales linearly in $n$ when cubic interpolation is used (and quadratic in $n$ when spectral interpolation is used), where  $n$ is the number of grid points in each velocity direction. Therefore, our methods offer significant computational savings over existing deterministic solvers in such cases. Numerical experiments on 2D and 3D benchmarks, including the BKW exact solution and anisotropic initial data, confirm the computational scaling, the expected order of accuracy  and verify the effectiveness of the conservation correction.
\end{abstract}

\textbf{Keywords:}
space-homogeneous Boltzmann equation; tensor train;   cross approximation; conservation

\section{Introduction}
 The Boltzmann equation is a fundamental equation in kinetic theory that describes
  the  evolution of the probability density function (\emph{pdf}) of a dilute gas through binary collisions. Its central
  computational challenge lies in the fast computation of non-linear and non-local collision operator, a five-dimensional integral that is notoriously expensive to evaluate.
 There have been two classes of numerical methods for solving the Boltzmann equation:  deterministic and stochastic. Stochastic methods are easy to implement but suffer from slow convergence and noise. In this paper, we will focus on the deterministic methods. Among the existing works, the fast Fourier spectral method \cite{pareschi2000numerical, bobylev1999fast, mouhot2006fast, gamba2017fast} is popular because it can reduce the computational cost to $\mathcal{O}(Mn^4 \log n)$ \cite{gamba2017fast} by exploiting the convolution structure of the collision operator, where $M$ is the number of angular quadrature points and $n$ is the number of grid points in each velocity direction. In \cite{gamba2009spectral} the authors proposed a method to enforce conservation of mass, momentum and energy based on Lagrange multiplier. The Petrov-Galerkin approach expresses the solution as combinations of spherical harmonics and Hermite/Laguerre polynomials \cite{gamba2018galerkin, wang2019approximation, kitzler2019polynomial}. The solution is defined on the whole $\mathbb{R}^{3}$ space, hence there is no need to choose a cut-off domain. In \cite{hanke2023representation}, the authors proposed a new approach that significantly reduces the memory cost and time cost, by taking advantage of the symmetry of collision operator. A different line of work is represented by discrete velocity methods \cite{mouhot2013convolutive, bobylev1995approximation,buet1996discrete}, which discretize the velocity variable directly and are valued for their robustness and their ability to preserve structural properties such as conservation, positivity, and entropy dissipation, at the price of higher computational cost and typically lower-order accuracy. However, the aforementioned work does not provide mechanism to deal with the curse of dimensionality inherent in the high-dimensional phase space of the Boltzmann equation.

To tackle the high dimensionality of the Boltzmann equation, in recent years, there have been tremendous development of low-rank methods, which explores the low-rank properties of the probability density function \cite{einkemmer2025review}. We mention \cite{einkemmer2021efficient,hu2022adaptive, dektor2025interpolatory, einkemmer2025asymptotic} which computes low-rank solutions to Boltzmann or BGK equations, among many others.
 However, most existing work focuses on the low-rank property between the velocity and the spatial variables, and the low-rank structure in the velocity space has not been explored for the Boltzmann collision operator due to the complex structure of the collision operator. On the other hand, we note  the low-rank property is well motivated because the equilibrium solution is low-rank in velocity space.

In this paper, we aim to develop a fast method to compute the collisional integral term exploring the low-rank property of the solution in the velocity space. The work is motivated by the recent seminal contribution in the analysis of Boltzmann and Landau equations \cite{guillen2025landau, imbert2026monotonicity, guillen2025landau2}, where a lifting-projection framework is introduced for analyzing the well-posedness of those equations.

We present a numerical version of this procedure: the \emph{lifting-projection (LP) scheme}
  for the Boltzmann equation. The key idea is to lift the original 3D equation to
  the 6D linear Kac master equation, evolve for one time step, and project the solution back to marginal distribution. The second contribution is a \emph{  low-rank
  tensor method} for evaluating the collision operator which computes with the lifted solution in 6D. We use the tensor train (TT) format \cite{oseledets2011tensor}  and present a numerical scheme in TT format utilizing interpolation and cross approximation, enhanced with conservation correction. The total computational cost per time
  step is $\mathcal{O}(nMmR^2 r^2 X + nR^3X)$, where $m$ is the width of interpolation stencil which is associated to the order of accuracy of the scheme, $M$ is the number of spherical quadrature points, $R$ and $r$ are the max TT ranks of the lifted and original \emph{pdf}, and $X$ is the number of sweeps in  cross approximation. We observe numerically  that $X\approx R$. Therefore, we can see that the computational cost scales linearly with $n,$ and when the ranks $R, r$ are small, this represents a  dramatic computational saving to existing deterministic methods. The fastest existing methods generally \cite{gamba2017fast} scale as $O(M n^{4} \log n)$, or as $O(M n^{3} \log n)$ in special cases \cite{mouhot2006fast}.

  The rest of the paper is organized as follows. In section~\ref{sec:boltzmann}, we introduce
  the Boltzmann equation and the LP scheme. In section~\ref{sec:TT}, we review
  the TT format and the cross approximation algorithms. In
  section~\ref{sec:LR}, we present the low-rank tensor method to compute the
  Boltzmann collision operator. In section~\ref{sec:num}, numerical results are
  presented to demonstrate the performance of the proposed method. Finally, we
  conclude the paper in section~\ref{sec:conclusion}.

\section{The space-homogeneous Boltzmann equation and its lifted version} \label{sec:boltzmann}
In this section, we review the space-homogeneous Boltzmann equation and its properties. We then introduce a framework of evaluating the collision operator by a lifting-projection procedure, which was inspired by the theoretical breakthrough of Guillen and Silvestre \cite{guillen2025landau} on well-posedness of the Landau equation. The same lifting technique was later adopted in the study of Boltzmann equation \cite{imbert2026monotonicity}.  A detailed expository article on this topic in kinetic theory can be found in \cite{guillen2025landau2}.

\subsection{Boltzmann equation and its properties}
The Boltzmann equation  describes the behavior of a dilute gas of particles when
the only interactions taken into account are binary elastic collisions. The space-homogeneous Boltzmann equation models the distribution function $f(\mathbf{v},t)$, which represents the probability density of finding a particle with velocity $\mathbf{v}$ at time $t$. The general form of the 3D Boltzmann equation is given by
\begin{equation}\label{eq:boltzmann}
    \partial_{t}f(\mathbf{v},t) = q(f)(\mathbf{v},t)
\end{equation}
where $q(f)$ is the collision operator,  defined as
\begin{equation*}
    q(f)(\mathbf{v},t) = \int_{\mathbb{R}^{3}}\int_{\mathbb{S}^{2}}\mathcal{B}(\mathbf{v}-\mathbf{w},\mathbf{\zeta})(f(\mathbf{v}')f(\mathbf{w}')-f(\mathbf{v})f(\mathbf{w}))\mathrm{d}\zeta\mathrm{d}\mathbf{w}.    
\end{equation*}
In this formula, $\mathbf{v}$ and $\mathbf{w}$ are the pre-collision velocities of two particles, $\mathbf{v}'$ and $\mathbf{w}'$ are the post-collision velocities, and $\mathcal{B}(\mathbf{v}-\mathbf{w},\mathbf{\zeta})$ is the collision kernel, which depends on the relative velocity $\mathbf{v}-\mathbf{w}$ and the scattering direction $\mathbf{\zeta} \in \mathbb{S}^{2}$ (the unit sphere). These quantities are related by the following formulas:
\begin{equation*}
    \begin{cases}
        \mathbf{v}' =& \frac{1}{2}\left(\mathbf{v} + \mathbf{w} + |\mathbf{v} - \mathbf{w}| \mathbf{\zeta}\right)\\
            \mathbf{w}' =& \frac{1}{2}\left(\mathbf{v} + \mathbf{w} - |\mathbf{v} - \mathbf{w}| \mathbf{\zeta}\right).
    \end{cases}
\end{equation*}

In this work, for simplicity, we will focus on the variable hard sphere (VHS) model, where the collision kernel is given by
\begin{equation*}
    \mathcal{B}(\mathbf{v}-\mathbf{w},\mathbf{\zeta}) = \frac{1}{4\pi}|\mathbf{v}-\mathbf{w}|^{\gamma},   
\end{equation*}
with $\gamma$ being a parameter that characterizes the type of interaction between particles. For instance, $\gamma = 1$ corresponds to hard spheres, while $\gamma = 0$ corresponds to Maxwell molecules. Nevertheless, the numerical method we propose in this work can be  extended to other collision kernels.

\begin{remark}
When $d=2$, the Boltzmann equation can be written in a similar form. The collision kernel in 2D is given by
\begin{equation*}
    \mathcal{B}(\mathbf{v}-\mathbf{w},\mathbf{\zeta}) = \frac{1}{2\pi}|\mathbf{v}-\mathbf{w}|^{\gamma}.
\end{equation*} 
and the scattering direction $\mathbf{\zeta} \in \mathbb{S}^{1}$, which can be parameterized by an angle $\theta \in [0, 2\pi)$.
\end{remark}

The Boltzmann equation has several important properties. Since the collision operator  satisfies the following properties:
\begin{equation*}
    \int_{\mathbb{R}^{3}}q(f)\begin{bmatrix}
1 \\
\mathbf{v}\\
|\mathbf{v}|^{2}
\end{bmatrix}d\mathbf{v} = 0,
\end{equation*}
i.e., the equation must conserve mass, momentum and energy: for any $t > 0$, we have
\begin{equation*}
    \int_{\mathbb{R}^{3}}f(\mathbf{v},t)\phi(\mathbf{v})d\mathbf{v} = \int_{\mathbb{R}^{3}}f(\mathbf{v},0)\phi(\mathbf{v})d\mathbf{v}, \forall~ \phi(\mathbf{v}) \in \{1, \mathbf{v}, |\mathbf{v}|^{2} \}.
\end{equation*}
Moreover, testing Equation \eqref{eq:boltzmann} with $\log f$ and integrate over $\mathbb{R}^{3}$, we have
\begin{equation*}
    \begin{split}
        \frac{d}{dt} \int_{\mathbb{R}^{3}}f\log f d\mathbf{v} =& \int_{\mathbb{R}^{3}}q(f)\log f d\mathbf{v}\\
        =& \int_{\mathbb{R}^{3}}\int_{\mathbb{R}^{3}}\int_{\mathbb{S}^{2}}\mathcal{B}(\mathbf{v}-\mathbf{w},\mathbf{\zeta})(f(\mathbf{v}')f(\mathbf{w}')-f(\mathbf{v})f(\mathbf{w}))\log f(\mathbf{v})d\zeta d\mathbf{w}d\mathbf{v}\\
        =& \frac{1}{4}\int_{\mathbb{R}^{3}}\int_{\mathbb{R}^{3}}\int_{\mathbb{S}^{2}}\mathcal{B}(\mathbf{v}-\mathbf{w},\mathbf{\zeta})(f(\mathbf{v}')f(\mathbf{w}')-f(\mathbf{v})f(\mathbf{w}))\log \frac{f(\mathbf{v})f(\mathbf{w})}{f(\mathbf{v}')f(\mathbf{w}')}d\zeta d\mathbf{w}d\mathbf{v}\\
        \leq& 0.
    \end{split} 
\end{equation*}
Hence the entropy defined as $H(f) \coloneqq \int_{\mathbb{R}^{3}}f\log f d\mathbf{v}$ is non-increasing in time.

\subsection{The scheme based on lifting-projection procedure} 
In \cite{guillen2025landau}, the authors resolved a long-standing open problem: the existence of global smooth solutions to the space-homogeneous Landau equation with Coulomb interaction. The key of the proof is to show decay of Fisher information, which was achieved through a lifting-projection procedure. The same lifting-projection technique was later extended to Boltzmann equation with very soft potentials \cite{imbert2026monotonicity}. Despite its huge impact on  analysis, we have not seen existing work that adopts this technique in numerical methods. This paper is an effort in this direction.

To illustrate the idea of the lifting-projection procedure, we first introduce the concept of tangent flow as a generalization of the forward Euler scheme. 

\begin{defi}[tangent flow]
    Let $(t, f(t)) \in \mathbb{R}^{+} \times V$ be the flow given by $\partial_{t} f = qf$ with initial condition $f^{(0)}$. 
    If there is another flow $(t, \widetilde{f}(t)) \in \mathbb{R}^{+} \times V$ such that
    \begin{equation*}
        \begin{split}
            f(\cdot, t_{0}) &= \widetilde{f}(\cdot, t_{0}),\\
            \partial_{t} f(\cdot, t_{0}) &= \partial_{t}\widetilde{f}(\cdot, t_{0}),\\
        \end{split}
    \end{equation*}
    then $(t,\widetilde{f}(t))$ is a tangent flow of $q$ at $(t_{0},f(t_{0}))$.
\end{defi}

\begin{remark}
Note that in particular, the forward Euler flow which is defined as $\widetilde{f}(t) = f(t_{0}) + (t-t_{0})q\left(f(t_{0})\right)$ is a tangent flow of $q$ at $(t_{0}, f(t_{0}))$. 
\end{remark}

For a function $F = F(\mathbf{v},\mathbf{w})$ defined on $\mathbb{R}^{3} \times \mathbb{R}^{3}$,  we define the lifted collision operator $Q$ by
\begin{equation*}
    Q(F)(\mathbf{v},\mathbf{w}) \coloneqq \int_{\mathbb{S}^{2}} \mathcal{B}(\mathbf{v}-\mathbf{w},\mathbf{\zeta}) \left(F(\mathbf{v}',\mathbf{w}') - F(\mathbf{v},\mathbf{w})\right)d\mathbf{\zeta}.
\end{equation*}

Let the projection operator $\Pi$ be defined as
\begin{equation*}
    \Pi F(\mathbf{v}) \coloneqq \int_{\mathbb{R}^{3}} F(\mathbf{v},\mathbf{w}) d\mathbf{w},  
\end{equation*}
and denote the tensor product of two functions as
\begin{equation*}
    (f\otimes g) (\mathbf{v}, \mathbf{w})\coloneqq f(\mathbf{v})g(\mathbf{w}),
\end{equation*}
then the Boltzmann collision operator $q(f)$ satisfies:
\begin{equation*}
    q(f) = \Pi Q (f\otimes f).
\end{equation*}
On the time interval $(t_{0}, t_{0}+T)$, suppose that $F(\mathbf{v},\mathbf{w},t)$ solve the lifted equation (a linear Kac master equation)
\begin{equation}\label{eq:2-kac}
    \partial_{t} F = Q (F)
\end{equation}
with initial condition $F(\mathbf{v},\mathbf{w},t_{0})=f(\mathbf{v},t_{0})f(\mathbf{w},t_{0})$, then we have
\begin{equation*}
    \left.\partial_{t} \Pi F(\mathbf{v},t) \right|_{t=t_{0}}= \Pi Q(f_{0}\otimes f_{0}) = \int_{\mathbb{R}^{3}} Q(f_{0}\otimes f_{0})d\mathbf{w}  = q(f_{0}),
\end{equation*}
which means the projection of the lifted solution $(t, \Pi F (\cdot, t))$ is a tangent flow of collision operator $q$ at $(t_{0}, f(\cdot, t_{0}))$. We call it the lifting-projection(LP) flow associated with $q$.

\begin{remark}\label{remark:aver}
    Define $\mathbf{u} = \mathbf{v} - \mathbf{w}$, $\mathbf{z} = \mathbf{v} + \mathbf{w}$, and $\mathbf{\sigma} = \mathbf{u}/|\mathbf{u}|$. Let $F = F(\lvert \mathbf{u} \rvert, \sigma; \mathbf{z},t)$, and define the spherical averaging operator as follows:
    \begin{equation}\label{eq:sph_aver}
        \mathcal{A} F(\lvert \mathbf{u} \rvert, \sigma;\mathbf{z}) \coloneqq \frac{1}{4\pi}\int_{\mathbb{S}^{2}}  F(\lvert \mathbf{u} \rvert, \zeta;\mathbf{z}) d\zeta,
    \end{equation}
    then the lifted Boltzmann operator asocciated to VHS model reads
    \begin{equation}\label{eq:lifted_VHS}
        QF = \lvert \mathbf{u} \rvert^{\gamma}\left(\mathcal{A}-I \right)F.
    \end{equation}
\end{remark}

Given a discrete solution $f_{h}^{(n)}$, we propose the following discrete LP scheme to obtain $f^{(n+1)}_{h}$.

\begin{enumerate}
    \item Solve the lifted Boltzmann equation in one time step
    \begin{equation}\label{eq:discrete_lifted}
    F^{(s+1)}_{h} = f^{(s)}_{h}\otimes f^{(s)}_{h} + \sqint_{s\Delta t}^{(s+1)\Delta t } Q_{h}F_{h}(\tau) d\tau,
    \end{equation}
    where $Q_{h}$ is the discretized version of lifted Boltzmann operator $Q$, $F_{h}(\tau)$ solves the equation $\partial_{t} F_{h} = Q_{h} F_{h}$ on $(s \Delta t, (s+1)\Delta t)$ with initial condition $f^{(s)}_{h}\otimes f^{(s)}_{h}$, and $\sqint$ represents the numerical integration in time. In this paper, we only consider the forward Euler time integrator, i.e. $\sqint_{s\Delta t}^{(s+1)\Delta t } Q_{h}F_{h}(\tau) d\tau = \Delta t Q_{h}F_{h}(s \Delta t)$
    \item Apply projection
    \begin{equation*}
        f^{(s+1)}_{h} =\Pi F^{(s+1)}_{h}.
    \end{equation*}
\end{enumerate}

\subsection{Time integration of the Boltzmann equation}
As have been shown in \cite{huang2025consecutive}, the  LP scheme introduced above is no more than first order accurate in time. To achieve higher order accuracy, we can adopt for example a high-order Runge-Kutta scheme. In this paper, we use the 4th order Runge-Kutta as follows:
\begin{equation}\label{eq:RK4}
    \begin{aligned}
        q_1 &= q\left(f_h^{(s)}\right),\\
        q_2 &= q\left(f_h^{(s)} + \tfrac{\Delta t}{2}\,q_1\right),\\
        q_3 &= q\left(f_h^{(s)} + \tfrac{\Delta t}{2}\,q_2\right),\\
        q_4 &= q\left(f_h^{(s)} + \Delta t\,q_3\right),\\
        f_h^{(s+1)} &= f_h^{(s)} + \frac{\Delta t}{6}\bigl(q_1 + 2q_2 + 2q_3 + q_4\bigr),
    \end{aligned}                        
\end{equation}
where
\begin{equation}\label{eq:subtract}
    q(g) = \Pi Q_h (g\otimes g)= \frac{1}{\Delta t} \left(\Pi (I+\Delta t Q_h)(g\otimes g) - g\right),
\end{equation}
which means we first solve the lifted equation \eqref{eq:discrete_lifted} by forward Euler method, project back to 3D, and then subtract the original function to obtain the time differencing. %The reason will be elaborated in \ref{subsec:alter_implement}

\section{Review of the TT format} \label{sec:TT}
Although the lifted Boltzmann equation is a linear equation, its solution is a 6D function, which makes it infeasible to store and compute with the full tensor representation when the grid size is large \cite{ballard2025tensor}. To overcome the curse of dimensionality, we use the tensor train (TT) format \cite{oseledets2011tensor} to represent the solution.

In this section, we review the TT format and the  cross approximation algorithm which will be used to solve the lifted equation. The TT format is a low-rank representation of high-dimensional tensors, which allows us to efficiently store and compute with high-dimensional data. The cross approximation algorithms are used to compute the TT representation without ever computing or storing the full tensor. Instead, we only need to sample   a small number of carefully chosen entries, which is essential for the efficiency of our low-rank Boltzmann solver.

\subsection{Definition and operations}
A tensor $A \in \mathbb{R}^{n_{1} \times n_{2} \times \cdots \times n_{d}}$ is said to be in TT  format \cite{oseledets2011tensor} if its element can be represented as: 
\begin{equation}\label{eq:TT_format}
    A(i_{1}, i_{2}, \cdots, i_{d}) = A^{(1)}_{\alpha_{0},\alpha_{1}}(i_{1})A^{(2)}_{\alpha_{1},\alpha_{2}}(i_{2})\cdots A^{(d)}_{\alpha_{d-1},\alpha_{d}}(i_{d})
\end{equation}
where $A^{(\mu)} \in \mathbb{R}^{r_{\mu-1} \times n_{\mu} \times r_{\mu}}$ are the cores, and $r_{\mu}$ are the ranks of the TT representation. When the ranks are low, the TT format offers tremendous saving for high-dimensional problems by reducing the storage cost from $O(n^d)$ to $O(dnr^2)$  where $n$ is the maximal dimension of the tensor and $r$ is the maximal rank of the cores.
We briefly review some of the basic operations of the TT format, while   refering the readers to \cite{oseledets2011tensor} for more details.

\subsubsection*{Summation}
The cores of the sum $C = A + B$ are given by
\begin{equation*}
    C^{(\mu)}(i_{\mu}) = \begin{bmatrix}
        A^{(\mu)}(i_{\mu})  & 0 \\
        0 & B^{(\mu)}(i_{\mu}) 
    \end{bmatrix},~\mu = 2, \cdots, d-1,
\end{equation*}
and 
\begin{equation*}
    C^{(1)}(i_{1}) = \begin{bmatrix}
        A^{(1)}(i_{1})  & B^{(1)}(i_{1})
    \end{bmatrix},~C^{(d)}(i_{d}) = \begin{bmatrix}
        A^{(d)}(i_{d})  \\
        B^{(d)}(i_{d})
    \end{bmatrix}.
\end{equation*}
As a result, the ranks of the sum are given by $r_{\mu}^{C} = r_{\mu}^{A} + r_{\mu}^{B},~\mu = 2, \cdots, d-1$.

\subsubsection*{Contraction}
The advantage of the TT format is that many operations can be implemented in a complexity proportional to $n$ instead of $n^{d}$. For example, consider the following one-dimensional contraction:
\begin{equation*}
    C(i_{1}, \cdots, i_{\mu-1}, i_{\mu+1}, \cdots, i_{d}) = \sum\limits_{i_{\mu}=1}^{n_{\mu}} A(i_{1}, i_{2}, \cdots, i_{d})u(i_{\mu}).  
\end{equation*}
If $A$ is in TT format as in Equation \eqref{eq:TT_format}, then the cores of $C$ are given by
\begin{equation*}
    C^{(1)} = A^{(1)}, \cdots, C^{(\mu-1)} = A^{(\mu-1)} \left(\sum\limits_{i_{\mu}=1}^{n_{\mu}} A^{(\mu)}(i_{\mu})u(i_{\mu}) \right)\in \mathbb{R}^{r_{\mu-2} \times r_{\mu}},~C^{(j-1)} = A^{(j)},~j = \mu+1, \cdots, d.
\end{equation*}
Multidimensional contraction can be implemented in a similar way with complexity $\mathcal{O}(nr^{2})$.

\subsubsection*{Rounding}
After performing operations on tensors in TT format (e.g. summation), the TT ranks may increase. To control the ranks and maintain the efficiency of the representation, one needs to perform a rounding procedure, which can be seen as a generalization of truncated singular value decomposition (SVD). The rounding procedure allows us to approximate the tensor with a lower rank representation while controlling the approximation error.
 The rounding procedure \cite{oseledets2011tensor} is listed in Algorithm \ref{alg:rounding}. It was shown in \cite{oseledets2011tensor} that the rounding procedure can be performed with complexity $\mathcal{O}(dnr^{3}).$

\begin{algorithm}[H]
\begin{algorithmic}[1]  
\Require {Tensor $A$ in TT format: $A(i_{1}, \cdots, i_{d}) = A^{(1)}_{\alpha_{0},\alpha_{1}}(i_{1})A^{(2)}_{\alpha_{1},\alpha_{2}}(i_{2}) \cdots A^{(d)}_{\alpha_{d-1},\alpha_{d}}(i_{d})$, tolerance $\varepsilon$.}
\Ensure {Lower rank tensor $B$ in TT format: $B(i_{1}, \cdots, i_{d}) = B^{(1)}_{\alpha_{0},\alpha_{1}}(i_{1})B^{(2)}_{\alpha_{1},\alpha_{2}}(i_{2}) \cdots B^{(d)}_{\alpha_{d-1},\alpha_{d}}(i_{d})$ such that}
\begin{equation*}
    \left\lVert A - B \right\rVert_{F} \leq \varepsilon \left\lVert A \right\rVert_{F}.
\end{equation*} 
\State{Compute truncation parameter $\delta = \frac{\varepsilon}{\sqrt{d-1}}\left\lVert A \right\rVert_{F}$.}
\State{Let $\widetilde{A}^{(d)}_{\alpha_{d-1}, \beta_{d}}(i_{d}) = A^{(d)}_{\alpha_{d-1}, \beta_{d}}(i_{d})$.}
\For{$\mu = d:-1:2$}   \Comment{Right-to-left orthonormalization.}
    \State{Perform QR decomposition on the 2D matrix $\widetilde{A}^{(\mu)}(\alpha_{\mu-1},\overline{i_{\mu}\beta_{\mu}})$:}
    \begin{equation*}
         R^{(\mu-1)}_{\alpha_{\mu-1},\alpha^{'}_{\mu-1}}Q^{(\mu)}_{\alpha^{'}_{\mu-1},\beta_{\mu}}(i_{\mu}) = \widetilde{A}^{(\mu)}_{\alpha_{\mu-1}, \beta_{\mu}}(i_{\mu}).
    \end{equation*}
    \State{Compute the new core $\widetilde{A}^{(\mu-1)}_{\alpha_{\mu-2},\beta_{\mu-1}} = A^{(\mu-1)}_{\alpha_{\mu-2}, \beta_{\mu-1}^{'}} R^{(\mu-1)}_{\beta_{\mu-1}^{'},\beta_{\mu-1}}$.}
\EndFor
\State{We have $A = \widetilde{A}^{(1)}_{\alpha_{0},\alpha_{1}} Q^{(2)}_{\alpha_{1},\alpha_{2}} \cdots Q^{(d)}_{\alpha_{d-1},\alpha_{d}}$.}
\State{Let $\widetilde{B}^{(1)}_{\beta_{0},\alpha_{1}}(i_{1}) = \widetilde{A}^{(1)}_{\beta_{0},\alpha_{1}}(i_{1})$.}
\For{$\mu=1:d-1$}  \Comment{Left-to-right truncation.}
    \State{Perform $\delta$-truncated SVD on the 2D matrix $\widetilde{B}^{(\mu)}(\overline{\gamma_{\mu-1} i_{\mu}}, \beta_{\mu})$:}
    \begin{equation*}
        B^{(\mu)}_{\gamma_{\mu-1},\gamma_{\mu}}(i_{\mu})S^{(\mu)}_{\gamma_{\mu},\gamma^{'}_{\mu}}V^{(\mu)}_{\gamma^{'}_{\mu},\beta_{\mu}} = \text{SVD}_{\delta}\left[\widetilde{B}^{(\mu)}_{\gamma_{\mu-1},\beta_{\mu}}(i_{\mu})\right].
    \end{equation*}
    \State{Compute the new core $\widetilde{B}^{(\mu+1)}_{\gamma_{\mu},\beta_{\mu+1}}(i_{\mu+1}) = S^{(\mu)}_{\gamma_{\mu},\gamma^{'}_{\mu}}V^{(\mu)}_{\gamma^{'}_{\mu},\alpha_{\mu}} Q^{(\mu+1)}_{\alpha_{\mu},\beta_{\mu+1}}(i_{\mu+1})$.}
\EndFor
\State{Return $ B (i_{1}, i_{2}, \cdots, i_{d}) = B^{(1)}_{\gamma_{0},\gamma_{1}}(i_{1})B^{(2)}_{\gamma_{1},\gamma_{2}}(i_{2}) \cdots \widetilde{B}^{(d)}_{\gamma_{d-1},\gamma_{d}}(i_{d})$.}
\end{algorithmic}
\caption{Rounding procedure for tensors in TT format. \cite{oseledets2011tensor}}
\label{alg:rounding}
\end{algorithm}

\subsection{Cross approximation for TT}

Cross approximation \cite{oseledets2010tt} is a class of algorithms for constructing a TT representation of a high-dimensional tensor accessing a few entries of the tensor, analogous to skeleton (CUR) decomposition for matrices. The key challenge is to select a small set of multi-indices at which to evaluate the tensor so that the resulting low-rank approximation is nearly optimal.

There are several cross approximation algorithms available in the literature \cite{dolgov2020parallel}, such as TT-Cross \cite{oseledets2010tt}, the DMRG greedy cross approximation algorithm \cite{SAVOSTYANOV2014217} etc. These algorithms differ in the way they select the indices for evaluation and update the cores, but they all aim to efficiently compute the TT representation with controlled ranks and approximation error. This is an active research area in multi-linear algebra. In this paper, we adopt the DMRG greedy cross method \cite{SAVOSTYANOV2014217}, which will be reviewed below.

The TT-cross approximation is given by:
\begin{equation}\label{eq:TTCA}
\begin{aligned}
&A(i_1,\dots,i_d)\\
\approx&\sum_{s_{1} \in \mathcal{I}^{\le 1},t_{1} \in \mathcal{I}^{>1}, \cdots,s_{d-1} \in \mathcal{I}^{\le d-1},t_{d-1} \in \mathcal{I}^{>d-1}}
A(i_1,t_1)
\left[ A(\mathcal{I}^{\le 1},\mathcal{I}^{>1}) \right]^{-1}(s_{1}, t_{1})\\
&A(s_1, i_2, t_{2})\left[ A(\mathcal{I}^{\le 2},\mathcal{I}^{>2}) \right]^{-1}(s_{2}, t_{2})
\cdots\cdots
A(s_{d-1}, i_d)
\\
=&\sum_{s_{k} \in \mathcal{I}^{\le k},t_{k} \in \mathcal{I}^{>k}, k = 0, 1, \cdots, d}
\left(\prod_{k=1}^{d}
A(s_{k-1}, i_k, t_{k})
\left[
A(\mathcal{I}^{\le k}, \mathcal{I}^{>k})
\right]^{-1}(s_{k}, t_{k})\right).
\end{aligned}
\end{equation}
Here $\mathcal{I}^{\le k}$ and $\mathcal{I}^{>k}$ denote the
positions of $r_k$ rows and columns in the $k$-th unfolding $A^{\{k\}}\coloneqq A(\overline{i_{1}\cdots i_{k}}, \overline{i_{k+1}\cdots i_{d}})$.
To unify the notation, we introduce the border sets
$\mathcal{I}^{\le 0} = \{1\}$ and
$\mathcal{I}^{> d} = \{1\}$.
The key step of the cross approximation algorithms is to select the indices $\mathcal{I}^{\le k}$ and $\mathcal{I}^{>k}$ for $k = 1, 2, \cdots, d-1$. The selection of these indices determines the accuracy and efficiency of the approximation.

In \cite{SAVOSTYANOV2014217}, the authors propose a greedy cross approximation algorithm that selects the indices based on the maximum-volume principle, which aims to maximize the determinant of the submatrix formed by the selected indices.
The algorithm considers interpolation sets $\mathcal{I}^{\le k}$ and $\mathcal{I}^{>k}$ that satisfy the following nestedness condition:
\begin{equation}\label{eq:nestedness}
\begin{split}
\overline{i_{1}\cdots i_{k}} \in \mathcal{I}^{\le k} &\Rightarrow \overline{i_{1}\cdots i_{k-1}} \in \mathcal{I}^{\le k-1}\\
\overline{i_{k+1}\cdots i_{d}} \in \mathcal{I}^{>k} &\Rightarrow \overline{i_{k+2}\cdots i_{d}} \in \mathcal{I}^{>k+1}
\end{split} 
\end{equation}
We describe the algorithm in Algorithm \ref{alg:greedy_cross} for completeness. The stopping criterion is that the relative $l_2$ error on a random test set is less than the given tolerance or the number of sweep reaches given threshold. This algorithm requires $\mathcal{O}(dnr^{2}X)$ evaluation of tensor elements and $\mathcal{O}(dnr^{3}X)$ additional operations, with $X$ being the number of sweeps. 

\begin{algorithm}[H]
\caption{Greedy restricted cross interpolation algorithm for TTs \cite{SAVOSTYANOV2014217}. }
\begin{algorithmic}[1]

\Require{Function to evaluate entries of a tensor 
$A(i_1,\dots,i_d)$, tolerance of error, max sweep.}

\Ensure{Cross interpolation (\ref{eq:TTCA}) with the nested interpolation sets (\ref{eq:nestedness}).}

\State{Choose $\mathcal{I}^{\le k}, \mathcal{I}^{>k}$, $k=1,\dots,d$, 
which satisfy (\ref{eq:nestedness}), and compute $\tilde{A}$ by (\ref{eq:TTCA}).}

\While{stopping criterion is not satisfied}

    \For{$k=1,\dots,d-1$} \Comment{Left-to-right half-sweep}
        \State{Apply cross interpolation \cite{oseledets2010cross} to the DMRG supercore matrix
        $[A(\mathcal{I}^{\le k-1}i_k, i_{k+1}\mathcal{I}^{>k+1})]$. Use interpolation sets $\mathcal{I}^{\le k}, \mathcal{I}^{>k}$ as the initial guess, and expand them to $\mathcal{J}^{\le k} \supset \mathcal{I}^{\le k}, \mathcal{J}^{>k} \supset \mathcal{I}^{>k}$, adding a few crosses:}
        \begin{equation*}
            A(\mathcal{I}^{\le k-1}i_k, i_{k+1}\mathcal{I}^{>k+1}) \approx \sum_{s_{k} \in \mathcal{J}^{\le k},t_{k} \in \mathcal{J}^{>k}}  
            A(\mathcal{I}^{\le k-1}i_k, t_{k})
            \left[A(\mathcal{J}^{\le k}, \mathcal{J}^{>k})\right]^{-1}(s_{k}, t_{k})
            A(s_{k}, i_{k+1}\mathcal{I}^{>k+1}).
        \end{equation*}

        \State{Substitute $\mathcal{I}^{\le k}, \mathcal{I}^{>k}$ by the expanded sets $\mathcal{J}^{\le k}, \mathcal{J}^{>k}$.}
    \EndFor

    \State{Perform right-to-left half-sweep in the same way}
\EndWhile
\end{algorithmic}
\label{alg:greedy_cross}
\end{algorithm}

\begin{remark}
    Let $X$ be the number of greedy DMRG sweeps, and let $r$ denote the final max TT rank of the approximated tensor. Empirically, each sweep only adds a small number of crosses, so in the tests below $X$ is comparable to the final rank scale, i.e. $X = \mathcal{O}(r)$.
\end{remark}

\section{Low-rank Boltzmann solver based on TT} \label{sec:LR}

With the TT format and the cross approximation algorithms, we are now ready to describe the fast scheme for the Boltzmann equation based on TT. For simplicity, we only describe the scheme for the 3D case. The 2D case is simpler, and the description  is skipped for brevity.

\subsection{Discretization}
As a standard practice, we need to choose a cut-off for the velocity domain.
For the 6D veclocity space, we choose $(\mathbf{v},\mathbf{w}) \in (-L, L)^{3} \times (-L, L)^{3}$ as the truncated domain. As have been illustrated in Remark \ref{remark:aver}, on each sphere $\mathcal{B}(\frac{z}{2},\lvert \mathbf{u} \rvert)$, the lifted Boltzmann drives the lifted solution toward spherical average through a reaction equation $F_{t} = \lvert \mathbf{u} \rvert^{\gamma} \left(\mathcal{A}F - F\right)$. When calculating the spherical average, we assume that the lifted solution $F$ is equal to zero outside of truncated domain.

We discretize the 3D truncated domain $(-L, L)^3\subset \mathbb{R}^{3}$ with a uniform grid of size $n$ in each dimension, and denote the grid points as $\{(v_{x}(i_{1}), v_{y}(i_{2}), v_{z}(i_{3}))\}$, where $i_{1}, i_{2}, i_{3} = 1, 2, \cdots, n$. In this paper, we assume the number of points in different directions to be the same, though the method is easily applicable with different number of points in different directions.  The grid points are given by
\begin{equation}\label{eq:gridpoint}
    \left(v_{x}(i_{1}), v_{y}(i_{2}), v_{z}(i_{3})\right) = \left(-L + (i_{1}-\frac{1}{2})h, -L + (i_{2}-\frac{1}{2})h, -L + (i_{3}-\frac{1}{2})h\right).
\end{equation}
where $h = \frac{2L}{n}$ is the grid size. Note that for simplicity we assume the numbers of nodes are the same in all directions. The discrete solution is denoted as $f_{h}(\mathbf{v})$, which is a tensor of size $n \times n \times n$:
\begin{equation*}
    A(i_{1}, i_{2}, i_{3}) = f_{h}\left(v_{x}(i_{1}), v_{y}(i_{2}), v_{z}(i_{3})\right).
\end{equation*}
The lifted solution is defined on the Cartesian product of the velocity grid with itself, and its discrete version is denoted as $F_{h}(\mathbf{v}, \mathbf{w})$, which is a tensor of size $n \times n \times n \times n \times n \times n$:
\begin{equation}\label{eq:xyordering}
    G(i_{1}, i_{2}, i_{3}, i_{4}, i_{5}, i_{6}) = F_{h}\left(v_{x}(i_{1}), v_{y}(i_{3}), v_{z}(i_{5}), w_{x}(i_{2}), w_{y}(i_{4}), w_{z}(i_{6})\right).
\end{equation}
We call this ordering the $xyz$-ordering, which  is the recommended ordering from our numerical experiment since it fully exploits the low-rank property between the velocity dimensions.

\begin{remark}
  The ordering of the dimensions in TT is critical for the performance of the solution in terms of ranks. For example, an alternative way to define the lifted solution (we call it $vw$-ordering) is:
    \begin{equation}
    \label{eq:vwordering}
        G(i_{1}, i_{2}, i_{3}, i_{4}, i_{5}, i_{6}) = F_{h}\left(v_{x}(i_{1}), v_{y}(i_{2}), v_{z}(i_{3}), w_{x}(i_{4}), w_{y}(i_{5}), w_{z}(i_{6})\right).
    \end{equation}
    We will show in the numerical experiments, (\ref{eq:xyordering}) leads to better low-rank structure of the lifted solution compared to \eqref{eq:vwordering}.
\end{remark}

When solving the VHS model, the key step is the computation of spherical averaging operator 
\begin{equation*}
    \mathcal{A} F(\lvert \mathbf{u} \rvert, \sigma;\mathbf{z}) \coloneqq \frac{1}{4\pi}\int_{\mathbb{S}^{2}}  F(\lvert \mathbf{u} \rvert, \zeta;\mathbf{z}) d\zeta,
\end{equation*}
where $\mathbf{u} = \mathbf{v} - \mathbf{w}$, $\mathbf{z} = \mathbf{v} + \mathbf{w}$, and $\mathbf{\sigma} = \mathbf{u}/|\mathbf{u}|$.

Denote the spherical quadrature points and weights as $\left\{(\mathbf{\zeta}_{i}, q_{i})\right\} \subset \mathbb{S}^{2} \times \mathbb{R}$, where $i = 1, 2, \cdots, M$, we define
\begin{equation*}
    A_{h} F_{h}(\mathbf{v},\mathbf{w}) = \sum\limits_{i=1}^{M} q_{i} F_{h}(\frac{\mathbf{z} + |\mathbf{u}|\mathbf{\zeta}_{i}}{2}, \frac{\mathbf{z} - |\mathbf{u}|\mathbf{\zeta}_{i}}{2}).
\end{equation*} 
Tensor product of one dimensional quadrature rules along the azimuthal and polar directions can be computationally expensive because the number of quadrature points grows quadratically with respect to the number of points in each direction. To address this issue, we use the Lebedev quadrature rule \cite{lebedev1999quadrature}, which is a family of quadrature rules for integration on the sphere that achieves high accuracy with a relatively small number of quadrature points.

Since quadrature points on the sphere does not necessarily lie on the discrete velocity grid, special treatment is needed. When it falls outside of the truncated domain $(-L, L)^{3}$, we simply use zero padding. For points inside the box $(-L,L)^{3}$, we use one of the following interpolations:
\begin{itemize}
    \item A 4-point stencil for interpolation in each dimension, which has order 4 accuracy. The interpolation weights are computed using the Lagrange polynomial basis. 
    \item Spectral interpolation with a global stencil using FFT.
\end{itemize}
More details will be described below.

\subsection{First-order in time scheme}
\begin{figure}[htbp]
    \centering
    \includegraphics[width=0.8\textwidth]{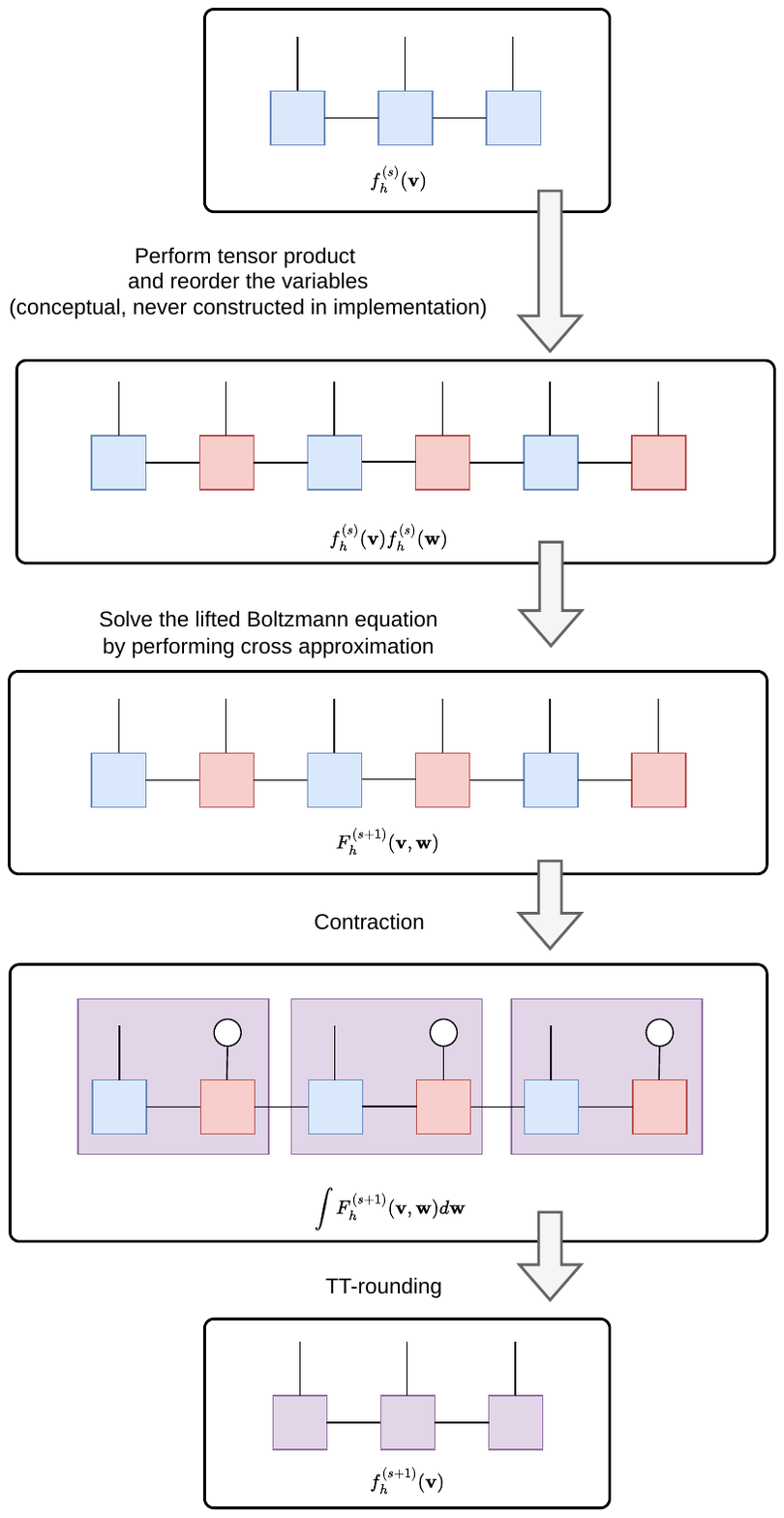}
    \caption{Flowchart of     Algorithm \ref{alg:updateTT}. }
    \label{fig:flowchart}
\end{figure}

Suppose we have the solution $f^{(s)}_{h}$ in TT format. To obtain the updated solution $f^{(s+1)}_{h}$ with first order accuracy in time, we take the following steps:  
\begin{enumerate}
    \item Compute the lifted solution $F_{h}^{(s+1)}$ in TT format using cross approximation (Algorithm \ref{alg:evalF}).
    \item Project the lifted solution $F_{h}^{(s+1)}$ to the original velocity space to obtain $f_{h}^{(s+1)}$ in TT format.
    \item Perform TT-rounding on $f_{h}^{(s+1)}$ to control the ranks of the solution. 
   % \item Perform conservation correction   (see Section \ref{sec:cons_weighted}).
\end{enumerate}
These steps are summarized in Algorithm \ref{alg:updateTT}, where a simple flowchart is plotted in Figure \ref{fig:flowchart}.

\begin{figure}[htbp]
    \centering
    \includegraphics[width=0.6\textwidth]{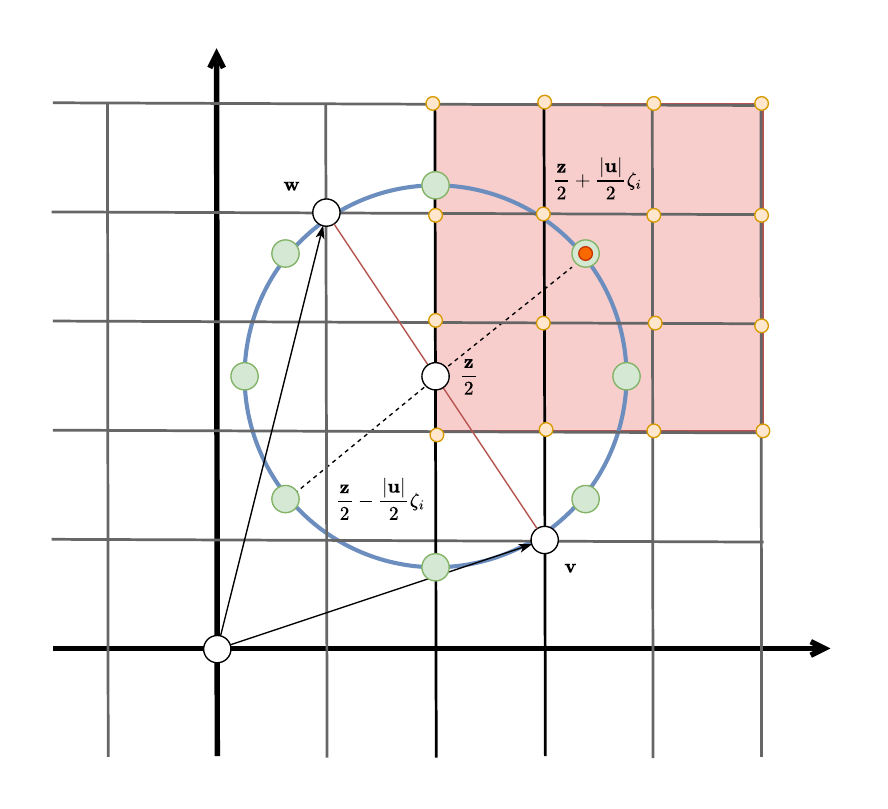}
    \caption{Visualization of evaluation of the 4D lifted solution associated with 2D Boltzmann equation. The green dots represent the quadrature points on the circle, and the orange dots represents the stencil for interpolation associated to one green dot using cubic interpolation.   }
    \label{fig:evalF}
\end{figure}

\begin{algorithm}[H]
\begin{algorithmic}[1]
\Require {The discrete solution $f^{(s)}_{h}$ in TT format $A(i_{1}, i_{2}, i_{3})=A^{(1)}_{\alpha_{0},\alpha_{1}}(i_{1})A^{(2)}_{\alpha_{1},\alpha_{2}}(i_{2})A^{(3)}_{\alpha_{2},\alpha_{3}}(i_{3})$}
\Ensure {The discrete solution $f^{(s+1)}_{h}$ in TT format $B(i_{1}, i_{2}, i_{3})=B^{(1)}_{\alpha_{0},\alpha_{1}}(i_{1})B^{(2)}_{\alpha_{1},\alpha_{2}}(i_{2})B^{(3)}_{\alpha_{2},\alpha_{3}}(i_{3})$.} 
\If{  spectral interpolation in Algorithm \ref{alg:evalF}}
    \State{Apply one-dimensional FFT to each TT core, obtaining $ \widehat{A}^{(1)}_{\alpha_{0},\alpha_{1}}$, $\widehat{A}^{(2)}_{\alpha_{1},\alpha_{2}}$ and $\widehat{A}^{(3)}_{\alpha_{2},\alpha_{3}}$.}
\EndIf
\State{Perform cross approximation using Algorithm \ref{alg:greedy_cross} by point evaluation from Algorithm \ref{alg:evalF} to obtain the TT form of the lifted solution:}
\begin{equation*}
    G(i_{1},\cdots,i_{6})=G^{(1)}_{\alpha_{0},\alpha_{1}}(i_{1})G^{(2)}_{\alpha_{1},\alpha_{2}}(i_{2})G^{(3)}_{\alpha_{2},\alpha_{3}}(i_{3})G^{(4)}_{\alpha_{3},\alpha_{4}}(i_{4})G^{(5)}_{\alpha_{4},\alpha_{5}}(i_{5})G^{(6)}_{\alpha_{5},\alpha_{6}}(i_{6})
\end{equation*}
\For{$ \ell = 1:3$}
        \State{Compute the marginal distribution in the $\ell$-th direction of $\mathbf{w}$.  $H^{(2\ell)}_{\alpha_{2\ell-1},\alpha_{2\ell}} = h \sum\limits_{i_{2\ell}=1}^{n}G^{(2\ell)}_{\alpha_{2\ell-1},\alpha_{2\ell}}(i_{2\ell})$}
        \State{ $B^{(\ell)}_{\alpha_{2\ell-2},\alpha_{2\ell}}(i_{2\ell-1})=\sum\limits_{\alpha_{2\ell-1}=1}^{r_{2\ell-1}}G^{(2\ell-1)}_{\alpha_{2\ell-2},\alpha_{2\ell-1}}(i_{2\ell-1})H^{(2\ell)}_{\alpha_{2\ell-1},\alpha_{2\ell}}$} 
\EndFor
\State{Perform TT-rounding on the tensor $B(i_{1}, i_{2}, i_{3})$.}

\end{algorithmic}
\caption{Low-rank scheme of first-order in time}
\label{alg:updateTT}
\end{algorithm}

Specifically, the cross approximation of the 6D lifted tensor is done by Greedy DMRG (Algorithm \ref{alg:greedy_cross}). The stopping criterion for the Greedy DMRG is that the relative $l_2$ error on a random test set is less than the given tolerance or the number of sweep reaches given threshold. In practice, we set the max number of sweep to be the same as grid size in each dimension. Moreover, the tolerance parameter is set to be on the order of the local truncation error, for example, $\frac{1}{2} \Delta t (h^{4} + (\Delta t)^{4})$ for cubic interpolation with RK4 scheme.    Point evaluations for the cross approximations are described in   Algorithm \ref{alg:evalF}, which is plotted in Figure \ref{fig:evalF} (for the 2D case with cubic interpolation for simplicity). The details of cubic and spectral interpolations are described in Algorithms \ref{alg:cubicInterp} and \ref{alg:spectralInterp}. For cubic interpolation (Algorithm \ref{alg:cubicInterp}), the computational cost scales as  $\mathcal{O}(Mr^{2})$.
For spectral interpolation (Algorithm \ref{alg:spectralInterp}), the computational cost scales as  $\mathcal{O}(Mr^{2}n)$. Finally, the purpose of the TT-rounding is to limit the rank of the TT. We suggest using the local truncation error as the rounding tolerance, for example, $\frac{1}{20} \Delta t (h^{4} + (\Delta t)^{4})$ for cubic interpolation with RK4 scheme.

\begin{algorithm}[H]
\begin{algorithmic}[1]
	\Require {Indices $(i_{1}, i_{2}, i_{3}, i_{4}, i_{5}, i_{6})$, the discrete solution $f^{(s)}_{h}$ in TT form $A^{n\times n \times n} = A^{(1)}_{\alpha_{0},\alpha_{1}}A^{(2)}_{\alpha_{1},\alpha_{2}}A^{(3)}_{\alpha_{2},\alpha_{3}}$ or its FFT $\widehat{A}^{n\times n \times n} = \widehat{A}^{(1)}_{\alpha_{0},\alpha_{1}}\widehat{A}^{(2)}_{\alpha_{1},\alpha_{2}}\widehat{A}^{(3)}_{\alpha_{2},\alpha_{3}}$, Lebedev nodes and weights $\left\{(\mathbf{\zeta}_{i}, q_{i})\right\} \subset \mathbb{S}^{2} \times \mathbb{R}$, $i = 1, 2, \cdots, M$, timestep $\Delta t$.  }
	\Ensure {Value of tensor element $G(i_{1}, i_{2}, i_{3}, i_{4}, i_{5}, i_{6})$ corresponding to the lifted solution $F_{h}^{(s+1)}$.} 
    \State{Evaluate $f_{h}(\mathbf{v}) = A(i_{1}, i_{3}, i_{5})$ and $f_{h}(\mathbf{w})= A(i_{2}, i_{4}, i_{6})$.}
    \State{Compute the variables $\mathbf{v}$ and $\mathbf{w}$ corresponding to the indices $(i_{1}, i_{2}, i_{3}, i_{4}, i_{5}, i_{6})$:}
    \begin{equation*}
        \begin{cases}
            \mathbf{v} &= \left(-L + (i_{1}-\frac{1}{2})h, -L + (i_{3}-\frac{1}{2})h, -L + (i_{5}-\frac{1}{2})h\right)\\
            \mathbf{w} &= \left(-L + (i_{2}-\frac{1}{2})h, -L + (i_{4}-\frac{1}{2})h, -L + (i_{6}-\frac{1}{2})h\right)
        \end{cases}
    \end{equation*}
    
    \State{Compute $\mathbf{z} = \mathbf{v} + \mathbf{w}$ and $\mathbf{u} = \mathbf{v} - \mathbf{w}$.}
	\State{Compute quadrature points $(\mathbf{v}_{i},\mathbf{w}_{i}) = \left(\frac{\mathbf{z} + |\mathbf{u}|\mathbf{\zeta}_{i}}{2}, \frac{\mathbf{z} - |\mathbf{u}|\mathbf{\zeta}_{i}}{2}\right)$, $i = 1, 2, \cdots, M$.}
    \If{using cubic interpolation}
        \State{Compute $f_{h}(\mathbf{v}_{i})$ and $f_{h}(\mathbf{w}_{i})$ with Algorithm \ref{alg:cubicInterp}.}
    \EndIf
    \If{using spectral interpolation}
        \State{Compute $f_{h}(\mathbf{v}_{i})$ and $f_{h}(\mathbf{w}_{i})$ with Algorithm \ref{alg:spectralInterp}.}
    \EndIf
    \State{$G(i_{1}, i_{2}, i_{3}, i_{4}, i_{5}, i_{6}) =  |\mathbf{u}|^{\gamma}\Delta t\sum\limits_{i=1}^{M}q_{i}f_{h}(\mathbf{v}_{i})f_{h}(\mathbf{w}_{i}) + (1-|\mathbf{u}|^{\gamma}\Delta t)f_{h}(\mathbf{v})f_{h}(\mathbf{w})$}
\end{algorithmic}
\caption{Pointwise evaluation of the lifted solution}
\label{alg:evalF}
\end{algorithm}

\begin{algorithm}[H]
\begin{algorithmic}[1]
\Require{TT cores:
  $A^{(\mu)}_{\alpha_{\mu}\alpha_{\mu+1}}$,\quad $\mu=1,2,3$; mean velocity $\frac{\mathbf{z}}{2}$ and radius of sphere $\frac{\lvert\mathbf{u}\rvert}{2}$;
  Lebedev nodes $\{\boldsymbol{\zeta}_i\}_{i=1}^M \subset \mathbb{S}^2$.
}
\Ensure{Values $f_{h}(\mathbf{v}_{i})$ where $\mathbf{v}_{i}= \frac{\mathbf{z} + |\mathbf{u}|\mathbf{\zeta}_{i}}{2}$.}

\For{$\mu = 1, 2, 3$}
  \State Form the sample points $x^{(\mu)}_i = \frac{1}{2}z^{(\mu)} + \frac{1}{2}\lvert\mathbf{u}\rvert\, \xi^{(\mu)}_i,\ i = 1,\cdots, M$.
  \State For each $\mu$ and $i$, compute the interpolation stencil $j_{k}^{(\mu)}(i)$ and the weights $W^{(\mu)}_{k}(i)$ with
    \begin{equation*}
        j_{k}^{(\mu)}(i) = \text{floor}\left(\frac{x^{(\mu)}_{i}-(-L+\frac{h}{2})}{h}\right) - 2 + k, \qquad k = 1,\cdots, 4,
    \end{equation*}
    and
    \begin{equation*}
    W^{(\mu)}(k,i)=\prod_{\substack{\ell=1\\ \ell\neq k}}^{4}
    \frac{x^{(\mu)}_i-v_{\mu}\!\left(j^{(\mu)}_{\ell}(i)\right)}
    {v_{\mu}\!\left(j^{(\mu)}_{k}(i)\right)-v_{\mu}\!\left(j^{(\mu)}_{\ell}(i)\right)},\qquad k=1,\ldots,4,
  \end{equation*}
  where $v_{\mu}(\cdot)$ is the grid point given by \eqref{eq:gridpoint}.
  
  \State Apply to the cores to get interpolated core slices: 
  \begin{equation*}
    \tilde{A}^{(\mu)}_{\alpha_{\mu}, \alpha_{\mu+1}}(i) =  \sum_{k = 1}^{4}A^{(\mu)}_{\alpha_{\mu}, \alpha_{\mu+1}}(j^{(\mu)}_{k}(i))\, W^{(\mu)}(k, i) \in \mathbb{R}^{r_{\mu} \times M \times r_{\mu+1} }.
  \end{equation*}
\EndFor

\State Evaluate at all $M$ Lebedev points:
\begin{equation*}
  f_h(\mathbf{v}_i) =
  \sum_{\alpha_1=1}^{r_1}\sum_{\alpha_2=1}^{r_2}
  \tilde{A}^{(1)}_{\alpha_{0}, \alpha_{1}} (i)\;
  \tilde{A}^{(2)}_{\alpha_{1}, \alpha_{2}} (i)\;
  \tilde{A}^{(3)}_{\alpha_{2}, \alpha_{3}} (i),
  \qquad i = 1,\ldots,M.
\end{equation*}

\end{algorithmic}
\caption{Cubic interpolation of a TT function at all Lebedev nodes (for 3D).}
\label{alg:cubicInterp}
\end{algorithm}

\begin{algorithm}[H]
\begin{algorithmic}[1]
\Require{
  FFT of each TT core along its grid axis:
  $\hat{A}^{(\mu)}_{\alpha_{\mu}\alpha_{\mu+1}}$,\quad $\mu=1,2,3$; mean velocity $\frac{\mathbf{z}}{2}$ and radius of sphere $\frac{\lvert\mathbf{u}\rvert}{2}$;
  Lebedev nodes $\{\boldsymbol{\zeta}_i\}_{i=1}^M \subset \mathbb{S}^2$. 
}
\Ensure{Values $f_{h}(\mathbf{v}_{i})$ where $\mathbf{v}_{i}= \frac{\mathbf{z} + |\mathbf{u}|\mathbf{\zeta}_{i}}{2}$.}

\For{$\mu = 1, 2, 3$}
  \State Form the sample points $x^{(\mu)}_i = \frac{1}{2}z^{(\mu)} + \frac{1}{2}\lvert\mathbf{u}\rvert\, \xi^{(\mu)}_i,\ i = 1,\cdots, M$.
  \State Build spectral weight matrix $W^{(\mu)} \in \mathbb{C}^{n \times M}$. 
  \begin{equation*}
      W^{(\mu)}(k,i)=\frac{1}{n}\exp\!\left(\mathrm{i}k \frac{2\pi}{2L}\left(x^{(\mu)}_{i}-(-L+\frac{h}{2})\right)\right),\qquad k=[n/2]+1-n,\ldots,[n/2],\quad i=1,\ldots,M.
  \end{equation*}
  \State Apply to FFT cores to get interpolated core slices:
  \begin{equation*}
    \tilde{A}^{(\mu)}_{\alpha_{\mu}, \alpha_{\mu+1}}(i) = \mathrm{Re}\left( \sum_{k = [n/2]+1-n}^{[n/2]}\hat{A}^{(\mu)}_{\alpha_{\mu}, \alpha_{\mu+1}}(k)\, W^{(\mu)}(k, i)\right) \in \mathbb{R}^{r_{\mu} \times M \times r_{\mu+1} }.
  \end{equation*}
\EndFor

\State Evaluate at all $M$ Lebedev points:
\begin{equation*}
  f_h(\mathbf{v}_i) =
  \sum_{\alpha_1=1}^{r_1}\sum_{\alpha_2=1}^{r_2}
  \tilde{A}^{(1)}_{\alpha_{0}, \alpha_{1}} (i)\;
  \tilde{A}^{(2)}_{\alpha_{1}, \alpha_{2}} (i)\;
  \tilde{A}^{(3)}_{\alpha_{2}, \alpha_{3}} (i),
  \qquad i = 1,\ldots,M.
\end{equation*}

\end{algorithmic}
\caption{Spectral interpolation of a TT function at all Lebedev nodes (for 3D).}
\label{alg:spectralInterp}
\end{algorithm}

\subsection{Higher-order in time and conservation correction}

The extensions to higher-order in time, e.g. RK4 is described in Algorithm \ref{alg:rk4}. The basic idea is to embed the first-order algorithm with a Runge-Kutta scheme. Here the tolerance of the TT-rounding is set to be on the order of the local truncation error of the scheme. In practice we are using $\frac{1}{20} \Delta t (h^{4} + (\Delta t)^{4})$ for cubic interpolation with RK4 scheme.  For numerical tests with spectral interpolation, we  set the TT-rounding and cross approximation tolerances to a small enough constant that can never be reached, which is $10^{-10}$ in our numerical experiments.  

\begin{algorithm}[H]
\caption{Low-rank schemes with RK4 and conservation correction}
\label{alg:rk4}
\begin{algorithmic}[1]
\Require Initial state $f_{h}^{(s)}$ in TT form, step size $\Delta t$, the forward Euler mapping $\Phi(f_{h}^{(s)},\Delta t) = \Pi (I + Q \Delta t) (f_{h}^{(s)}\otimes f_{h}^{(s)})$
\Ensure Next state $f_{h}^{(s+1)}$ in TT form
\State Compute $q_{1} = q(f^{(s)})= \frac{\Phi(f_{h}^{(s)},\Delta t) - f_{h}^{(s)}}{\Delta t}$ 
\State Compute $f_{1} = f^{(s)}_{h} + \frac{1}{2}\Delta t q_{1}$ and perform TT-rounding.
\State Compute $q_{2} = q(f_{1})= \frac{\Phi(f_{1},\Delta t) - f_{1}}{\Delta t}$ 
\State Compute $f_{2} = f^{(s)}_{h}  + \frac{1}{2}\Delta t q_{2}$ and perform TT-rounding.
\State Compute $q_{3} = q(f_{2})= \frac{\Phi(f_{2},\Delta t) - f_{2}}{\Delta t}$ 
\State Compute $f_{3} = f^{(s)}_{h}  + \Delta t q_{3}$ and perform TT-rounding.
\State Compute $q_{4} = q(f_{3})= \frac{\Phi(f_{3},\Delta t) - f_{3}}{\Delta t}$ 
\State Compute $f^{(s+1)}_{h} = f^{(s)}_{h} + \dfrac{\Delta t}{6}(q_1 + 2q_2 + 2q_3 + q_4)$ and perform TT-rounding.
\State Compute the correction term $g_{h}$ in TT form
\State Let $f_{h}^{(s+1)} = f_{h}^{(s+1)}+g_{h}$
\end{algorithmic}
\end{algorithm}

In addition, we develop a conservation correction step that is friendly to the TT format.
Due to numerical discretization, domain truncation and cross approximation, the conservation of mass, momentum and energy is generally lost. We propose to use a Lagragian multiplier approach and   solve the following constrained optimization problem to compute the correction term $\mu g_{h}$ to be added to the computed solution $f_{h}$:
\begin{equation*}
    \min_{g_{h}} \lVert g_{h}\rVert_{L^{2}_{\mu}}^{2} \text{ s.t. } (g_{h}, \phi_{k})_{L^{2}_{\mu}} = \mathcal{M}_{k} - (f_{h}, \phi_{k})_{L^{2}}, \quad k = 1, 2, \cdots, 5
\end{equation*}
where $\{\phi_{k}\} = \{1, v_{x}, v_{y}, v_{z}, |\mathbf{v}|^{2}\}$ to enforce mass, momentum and energy conservation \cite{gamba2009spectral}. The weight $\mu$ is set to be $\mu = \exp (-\lvert \mathbf{v}\rvert^{2}/2\sigma^{2})$, in practice we choose $\sigma = 1$. (When $\mu \equiv 1$ one recovers the correction in \cite{gamba2009spectral}.) We choose this weight to   accommodate the zero boundary conditions on the boundary, and in practice, it shows better performance in accuracy.

Let $A^{n^{3} \times 1}$ be the vectorization of the computed solution $f_{h}$, $M^{5 \times 1}$ be the conserved quantity,  $W^{n^{3}\times n^{3}}$ be the matrix associated to the weight $\mu$, and $\Phi^{n^{3}\times 5}$ be the matrix whose columns are the vectorization of $\{\phi_{k}\}$, then the optimization problem can be rewritten as 
\begin{equation*}
    \min_{G} G^{T}WG \text{ s.t. }  \Phi^{T} (W G + A) = M.
\end{equation*}
Analogous to \cite{gamba2009spectral}, the constrained optimization problem can be solved explicitly as follows:
\begin{equation*}
    G = \Phi(\Phi^{T}W\Phi )^{-1}( M - \Phi^{T}A ).
\end{equation*}

In implementation,   this correction is performed in the TT format. Since $\{\phi_{k}\}$ are low-rank functions in the velocity space, $\Phi$ can be represented in TT format with low rank. Consequently, the correction term $WG$ can be computed with complexity $\mathcal{O}(dnr^{2})$ instead of $\mathcal{O}(n^{3})$ for the full tensor. %The details of the algorithm are presented in Algorithm \ref{alg:updateTT}.
In addition, we note the correction is of low-rank in TT format. This is because  $g_{h}$ is always a linear combination of $\{\phi_{k}\} = \{1, v_{x}, v_{y}, v_{z}, \mathbf{v}^{2}\}$.  If we define 
\begin{equation*}
    \begin{split}
        &g_{1}(v_{x},1) = 1,\ g_{1}(v_{x},2) = \kappa_{1}v_{x} + \kappa_{4}v_{x}^{2},\\
        &g_{2}(v_{y},1,1) = 1 + \kappa_{2}v_{y} + \kappa_{4}v_{y}^{2},\ g_{2}(v_{z},2,2) = 0, g_{2}(v_{y}, 1,2) = g_{2}(v_{y}, 2,1)=1,\\
        &g_{3}(v_{z},1) = 1,\ g_{3}(v_{z},2) = \kappa_{3}v_{z} + \kappa_{4}v_{z}^{2},
    \end{split}
\end{equation*}
then
\begin{equation*}
    1 + \kappa_{1}v_{x} + \kappa_{2}v_{y} + \kappa_{3}v_{z} + \kappa_{4} (v_{x}^{2} + v_{y}^{2} + v_{z}^{2})  = \sum_{\alpha_{1} = 1}^{2}\sum_{\alpha_{2} = 1}^{2} g_{1}(v_{x},\alpha_{1})g_{2}(v_{y},\alpha_{1}, \alpha_{2})g_{3}(v_{z},\alpha_{2}).
\end{equation*}
and it follows that
\begin{equation*}
    \begin{split}
        &\left[1 + \kappa_{1}v_{x} + \kappa_{2}v_{y} + \kappa_{3}v_{z} + \kappa_{4} (v_{x}^{2} + v_{y}^{2} + v_{z}^{2})\right]\mu_{1}(v_{x})\mu_{2}(v_{y})\mu_{3}(v_{z}) \\ =& \sum_{\alpha_{1} = 1}^{2}\sum_{\alpha_{2}=1}^{2} \mu_{1}g_{1}(v_{x},\alpha_{1})\mu_{2}g_{2}(v_{y},\alpha_{1}, \alpha_{2})\mu_{3}g_{3}(v_{z},\alpha_{2}).
    \end{split}
\end{equation*}
Therefore, the correction term $\mu g_{h}$ is always a TT with max rank equal to 2.

\subsection{Complexity analysis}
Suppose that the lifted 6D solution has TT-ranks $\{R_{1}, \cdots, R_{5}\}$ with maximum $R$, and the projected 3D solution has TT-ranks $\{r_{1},r_{2}\}$ with maximum $r$.  

The cost of TT cross approximation is dominant in Algorithm \ref{alg:updateTT}. The greedy cross approximation Algorithm \ref{alg:greedy_cross} takes $\mathcal{O}(nR^{2}X)$ evaluations of the tensor elements, and $\mathcal{O}(nR^{3}X)$ flops for other operations, with $X$ being the number of sweeps.  
It takes $\mathcal{O}(Mr^{2})$   interpolations for each query of   Algorithm \ref{alg:evalF}. Each cubic interpolation (Algorithm \ref{alg:cubicInterp}) takes $\mathcal{O}(1)$ flops and each spectral interpolation (Algorithm \ref{alg:spectralInterp}) takes $\mathcal{O}(n)$ flops. 
Therefore the total cost of one time step is in the order of $\mathcal{O}(nMmR^{2}r^{2}X + nR^{3}X)$, where $m$ is the width of interpolation stencil.

Finally, we would like to briefly comment on alternative implementations based on TT.   We have conducted numerical investigations into two alternative approaches, and below is our findings.

\begin{itemize}
    \item (Alternative 1) Perform cross approximation directly on the lifted operator $Q(g\otimes g)$ rather than $\left(I+ \Delta t Q\right)(g\otimes g)$.
The solver does not perform well. This is because when $g$ is almost in equilibrium, $Q(g \otimes g )$ is near zero, which is not computational robust for cross approximation. 

\item (Alternative 2) Directly compute the original Boltzmann equation \eqref{eq:boltzmann} with TT without lifting-projection. If we perform a TT cross approximation on the forward Euler scheme for \eqref{eq:boltzmann}, this will not take advantage of   fast operations of the lifted solution. For example, the cross approximation Algorithm \ref{alg:greedy_cross} needs $\mathcal{O}(nr^{2})$ point evaluations of tensor elements. 

For each 3D velocity $\mathbf{v}$, evaluation of $\Pi(I + Q \Delta t) (g\otimes g)(\mathbf{v})$ requires $\mathcal{O}(M n^{3})$ quadrature points since $\Pi Q$ includes integration on $\mathbb{R}^{3}\times \mathbb{S}^{2}$. If $g$ is given as a 3D TT with max rank $r$, then the total cost will be $\mathcal{O}(Mn^{4}r^{4})$, considering that the cubic interpolation of $g$ at off-grid velocity takes $\mathcal{O}(r^{2})$ flops. In most cases this is much more expensive than the approach we propose.
\end{itemize}

\section{Numerical results} \label{sec:num}

In this section, we conduct numerical simulations on 2D and 3D benchmark problems, focusing on the performance of the method in terms of accuracy, speed and low-rank properties. 
In what follows, we choose $L=6.4,$ i.e. using $(-6.4, 6.4)^{2}$ or $(-6.4, 6.4)^{3}$ as the truncated domain.

\subsection{Experiments on 2D Boltzmann equation}
We use the 2D case to study the effect of the ordering for the ranks and behaviors of the solution, as well as benchmarking on numerical accuracy.

\subsubsection{Ordering of variables for TT}

As has been shown in Section \ref{sec:LR}, the complexity of the low-rank method depends on the rank of TTs. Therefore, it is crucial to investigate how the ordering of variables affects the TT rank of the lifted solution. We consider:
\begin{enumerate}
    \item The $vw$-ordering $(v_{x}, v_{y}, w_{x}, w_{y})$ (as defined in \eqref{eq:vwordering}).
    \item The $xy$-ordering $(v_{x}, w_{x}, v_{y}, w_{y})$ (as defined in \eqref{eq:xyordering} for the 3D case).
\end{enumerate}
The orderings are tested on three different initial conditions with the Maxwell molecules where $\gamma = 0$:
\begin{enumerate}
    \item The BKW \cite{bobylev1975exact, krook1977exact} reference solution:
        \begin{equation}
        \label{eq:bkw2d}
            f_{0}(\mathbf{v}) = |\mathbf{v}|^{2}\exp(-|\mathbf{v}|^{2})
        \end{equation}
    \item Two Gaussians centered at $\mathbf{v}_{1} = (-2, 0)$ and $\mathbf{v}_{2} = (2, 0)$:
        \begin{equation}
        \label{eq:2gaussian2d}
            f_{0}(\mathbf{v})=\exp(-\frac{|\mathbf{v}-\mathbf{v}_{1}|^2}{2}) + \exp(-\frac{|\mathbf{v}-\mathbf{v}_{2}|^2}{2})
        \end{equation}
    \item Anisotropic Gaussian centered at the origin:
        \begin{equation}
        \label{eq:anigaussian2d}
            f_{0}(\mathbf{v}) = \exp(-\frac{(v_{x} + v_{y})^{2}}{8}-\frac{(v_{x} - v_{y})^{2}}{2})
        \end{equation}
\end{enumerate}
We note that the last example is more challenging for low-rank methods due to its intrinsic higher rank between $v_{x}$ and $v_{y}$.  

In each test, let the mesh size $h$ be $\frac{2L}{n} = \frac{2\times6.4}{32} = 0.4$. For purpose of benchamarking, we use full tensor, i.e. in Algorithm \ref{alg:updateTT}, we keep the same interpolation and quadrature rules, but replace the TT cross-approximation step with full-tensor construction. We compute the full 4D  tensor   associated to the lifted solution at the first time step $t = \Delta t = 0.1$, denoted by $G$, and perform TT-SVD to obtain the optimal low-rank approximation $G_{*}$ such that:
\begin{equation*}
    \left\lVert G - G_{*} \right\rVert_{F} \leq \varepsilon \left\lVert G \right\rVert_{F},
\end{equation*}  
where $\varepsilon$ is set to be $10^{-5}$.

The results are summarized in Table \ref{tab:ttrank}, where we report the TT ranks of the lifted solution $G$ with respect to various orderings and initial conditions. It is clear the advantage of the $xy$-ordering compared to $vw$-ordering for BKW and two Gaussian test cases where the true solution exhibits low-rank properties. For the anisotropic Gaussian case, the two methods of ordering show similar ranks. Since the main use case for low-rank methods are when the solution $f$ exhibits low-rank properties, the $xy$-(or $xyz$-ordering in 3D) is well motivated.

\begin{table}[H]
    \centering
    \begin{tabular}{cccc}
         &  BKW&  Two Gaussians& Anisotropic Gaussian\\
         $vw$-ordering&  (14, 92, 14)&  (18,135,16)& (16, 84, 16) \\
         $xy$-ordering&  (14, 22, 14)&  (18,16,16) & (16, 95, 16) \\
    \end{tabular}
    \caption{Rank $(r_{1}, r_{2}, r_{3})$ of the 4d tensor associated with the three test cases and the two way of ordering variables.}
    \label{tab:ttrank}
\end{table}

\subsubsection{Accuracy test}
To verify the order of accuracy, we present the results of the following test: Using maxwell molecule model, i.e. $\mathcal{B}(\mathbf{v}-\mathbf{w},\zeta) = \frac{1}{2\pi}$, we let the BKW initial condition
    \begin{equation*}
        f(\mathbf{v}, t) = \frac{1}{2\pi K(t)^{2}}\exp\left(-\frac{|\mathbf{v}|^{2}}{2K(t)}\right)\left( 2K(t)-1 + \frac{1-K(t)}{2K(t)}|\mathbf{v}|^{2}\right),
    \end{equation*}
    where $K(t) = 1 - 0.5\exp{(-t/8)}$, evolve from $t=0$ to $t=4$, and compute the relative $L^{2}$ error by comparing the computed solution with the reference solution given by the BKW formula. We use RK4 as time integrator with cubic interpolation. The run uses $M=64$ circular quadrature points, cross-approximation tolerance $\frac{1}{2}\Delta t(h^4+(\Delta t)^4)$, at most $n$ greedy sweeps, and TT-rounding tolerance $10^{-16}$.
    
The results are shown in Table \ref{tab:2d_accuracy_mm}. We can clearly observe 4th order convergence in space and time.
\begin{table}[H]
    \centering
    \begin{tabular}{c|c|c|c|c|c|c}
        $(h, \Delta t)$ & $t=1.2$& order & $t=2$& order & $t=4$& order\\
        $(0.4, 0.4)$ & 2.8e-3& -& 3.5e-3& - & 2.8e-3& -\\
        $(0.2, 0.2)$ & 2.0e-4& 3.8& 2.1e-4& 4.1& 1.7e-4& 4.0\\
        $(0.1, 0.1)$ & 1.2e-5& 4.1& 1.4e-5& 3.9& 1.1e-5& 3.9\\
    \end{tabular}
    \caption{Relative $L^2$ error for maxwell molecule model test case with cubic interpolation for 2D. We calculate the error by comparing the computed solution with the BKW reference solution.}
    \label{tab:2d_accuracy_mm}
\end{table}

\subsection{Experiments on 3D Boltzmann equation}
\subsubsection{Accuracy test on cubic interpolation}
First we test the RK4 method: Algorithm \ref{alg:rk4} with cubic interpolation: Algorithm \ref{alg:cubicInterp}. In these tests we use $M=74$ Lebedev points, cross-approximation tolerance $\frac{1}{2}\Delta t(h^{4}+(\Delta t)^4)$, a maximum number of greedy sweeps equal to the number of grid points in each velocity direction, and TT-rounding tolerance equal to one tenth of the cross-approximation tolerance.

 We expect the overall order of accuracy of the proposed method to be 4. To verify the order of accuracy, we present the results of two experiments:
\begin{itemize}
    \item  Maxwell molecule model, i.e. $\mathcal{B}(\mathbf{v}-\mathbf{w},\zeta) = \frac{1}{4\pi}$, we consider the BKW solution
    \begin{equation}
    \label{eq:bkw3d}
        f(\mathbf{v}, t) = \frac{1}{(2\pi K(t))^{3/2}}\exp\left(-\frac{|\mathbf{v}|^{2}}{2K(t)}\right)\left( \frac{5K(t)-3}{2K(t)} + \frac{1-K(t)}{2K(t)^{2}}|\mathbf{v}|^{2}\right),
    \end{equation}
    where $K(t) = 1 - \exp{(-(t+5.5)/6)}$, evolve from $t=0$ to $t=2$, and compute the relative $L^{2}$ error by comparing the computed solution with the reference solution given by the BKW formula. 
    \item Hard sphere model, i.e. $\mathcal{B}(\mathbf{v}-\mathbf{w},\zeta) = \frac{1}{4\pi} \left\vert\frac{\mathbf{v}-\mathbf{w}}{6.4}\right\vert$, we let the  two Gaussians initial condition 
    \begin{equation} \label{eq:twoG3d}
        f_{0}(\mathbf{v})=\exp(-\frac{|\mathbf{v}-\mathbf{v}_{1}|^2}{2}) + \exp(-\frac{|\mathbf{v}-\mathbf{v}_{2}|^2}{2})
    \end{equation}
    where $\mathbf{v}_{1} = (2, 0, 0)$ and $\mathbf{v}_{2} = (0, -2, 0)$, evolve from $t=0$ to $t=2$. Since there is no exact reference solution, we only perform  self-convergence test. 
\end{itemize}

The results are shown in Tables \ref{tab:accuracy_cubic_mm} and \ref{tab:accuracy_cubic_hs}. In both cases, we can clearly observe 4th order convergence upon mesh refinement.  
\begin{table}[H]
    \centering
    \begin{tabular}{c|c|c|c|c|c|c}
        $(h, \Delta t)$ & $t=0.8$ & order & $t=1.2$ & order & $t=2$ & order\\
        $(0.4, 0.4)$ & 2.0e-3 & -& 2.5e-3 & - & 2.0e-3 & -\\
        $(0.2, 0.2)$ & 1.0e-4 & 4.3& 9.1e-5 & 4.8& 8.9e-5 & 4.5\\
        $(0.1, 0.1)$ & 5.0e-6 & 4.3& 5.8e-6 & 4.0 & 6.0e-6 & 3.9\\
    \end{tabular}
    \caption{Relative $L^2$ error for maxwell molecule model test case with cubic interpolation. We calculate the error by comparing the computed solution with the BKW reference solution.} \label{tab:accuracy_cubic_mm}
    \begin{tabular}{c|c|c|c|c|c|c}
        $(h, \Delta t)$ & $t=0.8$ & order & $t=1.2$ & order & $t=2$ & order\\
        $(0.4, 0.4)$ & 9.8e-3 & -& 1.7e-2 & - & 2.6e-2 & -\\
        $(0.2, 0.2)$ & 3.5e-4 & 4.8 & 4.2e-4 & 5.3& 2.0e-3 & 3.7
    \end{tabular}
    \caption{Relative $L^2$ error for hard sphere model test case with cubic interpolation. We calculate the error by comparing the computed solution with the solution obtained by $(h,\Delta t ) = (0.1, 0.1)$. }\label{tab:accuracy_cubic_hs}
\end{table}

\begin{figure}[H]
    \centering
    \begin{subfigure}[b]{0.45\textwidth}
        \centering
        \includegraphics[width=\textwidth]{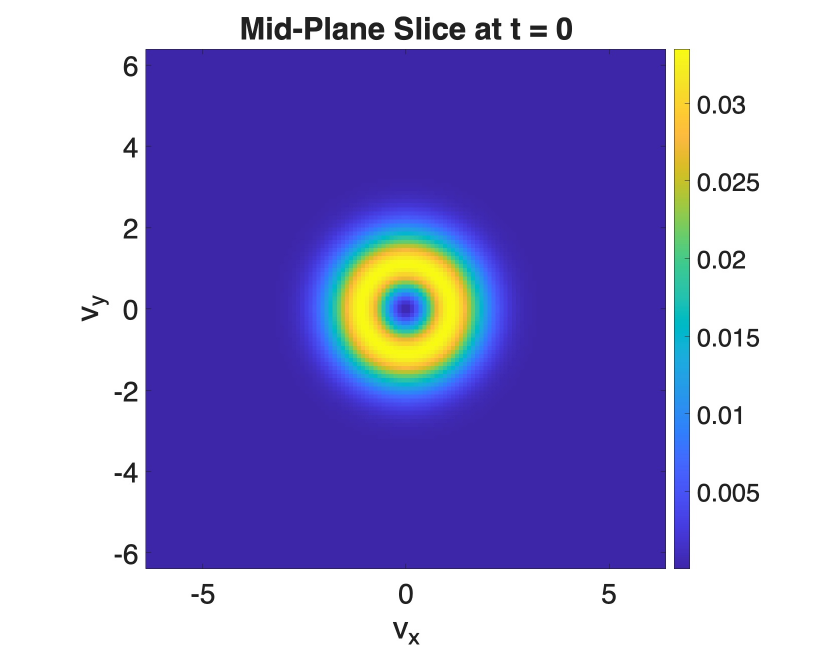}
        \caption{Initial state at $t = 0$.}
    \end{subfigure}
    \hfill
    \begin{subfigure}[b]{0.45\textwidth}
        \centering
        \includegraphics[width=\textwidth]{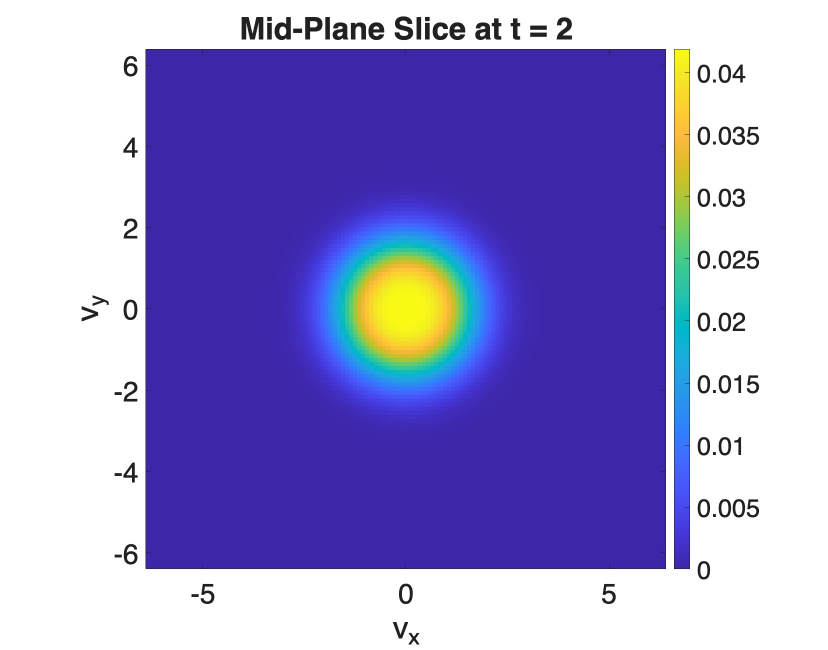}
        \caption{Final state at $t = 2$.}
    \end{subfigure}
    \caption{Slice of the computed solution for maxwell molecule model with BKW initial condition. The solution is sliced at $v_{z} = 0$.}\label{fig:plot_mm}
    \begin{subfigure}[b]{0.45\textwidth}
        \centering
        \includegraphics[width=\textwidth]{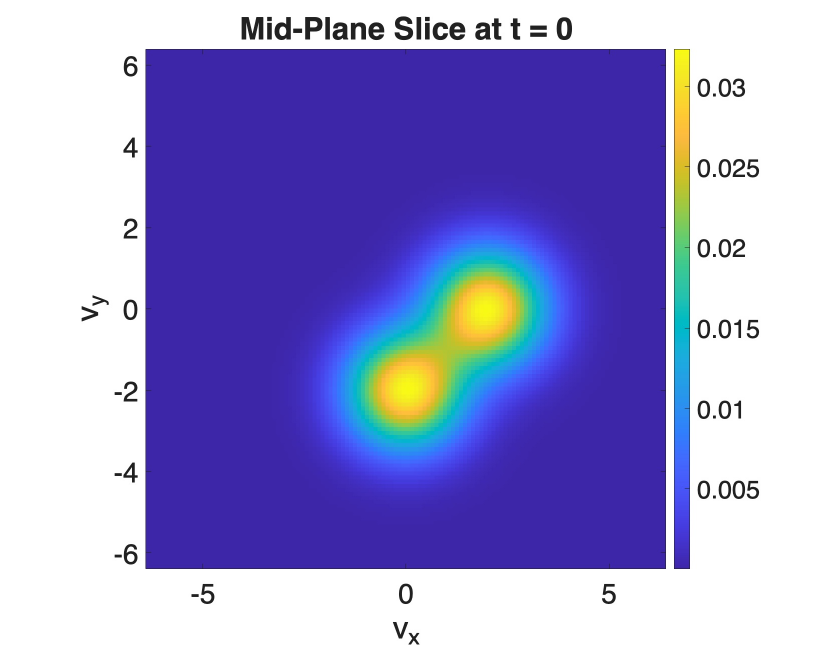}
        \caption{Initial state at $t = 0$.}
    \end{subfigure}
    \hfill
    \begin{subfigure}[b]{0.45\textwidth}
        \centering
        \includegraphics[width=\textwidth]{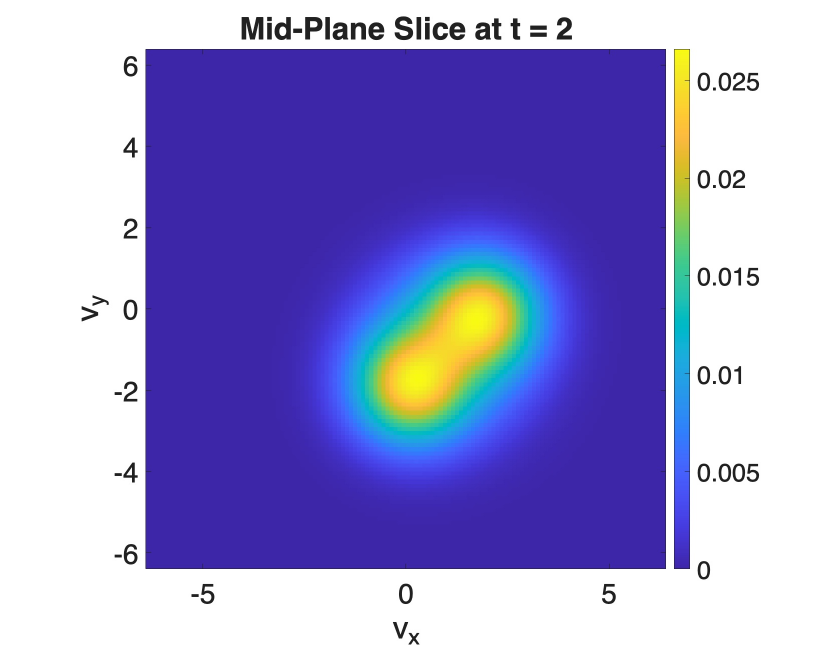}
        \caption{Final state at $t = 2$.}
    \end{subfigure}
    \caption{Slice of the computed solution for hard sphere model with two Gaussians initial condition. The solution is sliced at $v_{z} = 0$.}\label{fig:plot_hs}
\end{figure}

In addition, the initial and final states are plotted in Figures \ref{fig:plot_mm} and \ref{fig:plot_hs}, which benchmark well with known solutions in the literature. 

\subsubsection{Accuracy test on spectral interpolation}
Next we test the RK4 method: Algorithm \ref{alg:rk4} with spectral interpolation: Algorithm \ref{alg:spectralInterp}.

We expect the method to have spectral accuracy in velocity and 4th order accuracy in time. We solve the same test cases as in \eqref{eq:bkw3d} and \eqref{eq:twoG3d}. For the spectral runs we use $M=74$ Lebedev points, cross-approximation tolerance $10^{-10}$, at most $n$ greedy sweeps, and TT-rounding tolerance $10^{-10}$. 

The results of Maxwell molecules with BKW initial condition are shown in Table \ref{tab:accuracy_mm_spectral}. As we can see, extremely accurate numerical solution can be obtained with $n=32,$  reflecting spectral accuracy.  The error saturates at the level of $10^{-8}$ due to the prescribed truncation tolerance and possibly the temporal discretization error.

The results of hard spheres with ``two-Gaussians" initial condition are shown in Table \ref{tab:accuracy_hs_spectral}. In this case, the solution is not strictly low-rank. We do not observe spectral accuracy because  the cross approximation algorithm reaches max number of sweep and fails to meet the tolerance.

\begin{table}[H]
    \centering
    \begin{tabular}{c|c|c}
        $(n, \Delta t)$ & $t=0.5$ & $t=1.0$ \\
        $(16, 0.1)$ & 7.9e-4 & 1.0e-3 \\
        $(32, 0.1)$ & 1.5e-8 & 2.0e-8 \\
        $(64, 0.1)$ & 1.5e-8 & 1.9e-8 \\
    \end{tabular}
    \caption{Relative $L^2$ error for maxwell molecule model test case with spectral interpolation. We calculate the error by comparing the computed solution with the BKW reference solution.}\label{tab:accuracy_mm_spectral}
    \begin{tabular}{c|c|c}
        $(n, \Delta t)$ & $t=0.5$ & $t=1.0$ \\
        $(16, 0.1)$ & 2.3e-2 & 2.1e-2 \\
        $(32, 0.1)$ & 2.4e-3 & 2.2e-3 
    \end{tabular}
    \caption{Relative $L^2$ error for hard sphere model test case with spectral interpolation. We calculate the error by comparing the computed solution with the solution obtained by $(n,\Delta t ) = (64,0.1)$.}\label{tab:accuracy_hs_spectral}
\end{table}

\subsubsection{Complexity of cubic interpolation}\label{subsubsec:cubic_complexity}

Recall that with cubic interpolation, the complexity of the proposed method is $\mathcal{O}(nMR^{2}r^{2}X + n R^3 X)$, where $n$ is the number of grid points in each direction, $M$ is the number of quadrature points on the sphere, $R$ is the max TT rank of the 6d lifted solution $F$, $r$ is the TT-rank of the 3D projected solution $f$, and $X$ is the number of sweeps. In what follows we present the computational complexity  with two test cases, \eqref{eq:bkw3d} and \eqref{eq:twoG3d}. The timing runs use $M=74$ Lebedev points, a maximum number of greedy sweeps equal to the number of grid points per velocity direction, and TT-rounding tolerance equal to one tenth of the cross-approximation tolerance. 

\begin{figure}[H]
    \centering
    \begin{subfigure}[b]{0.45\textwidth}
        \centering
        \includegraphics[width=\textwidth]{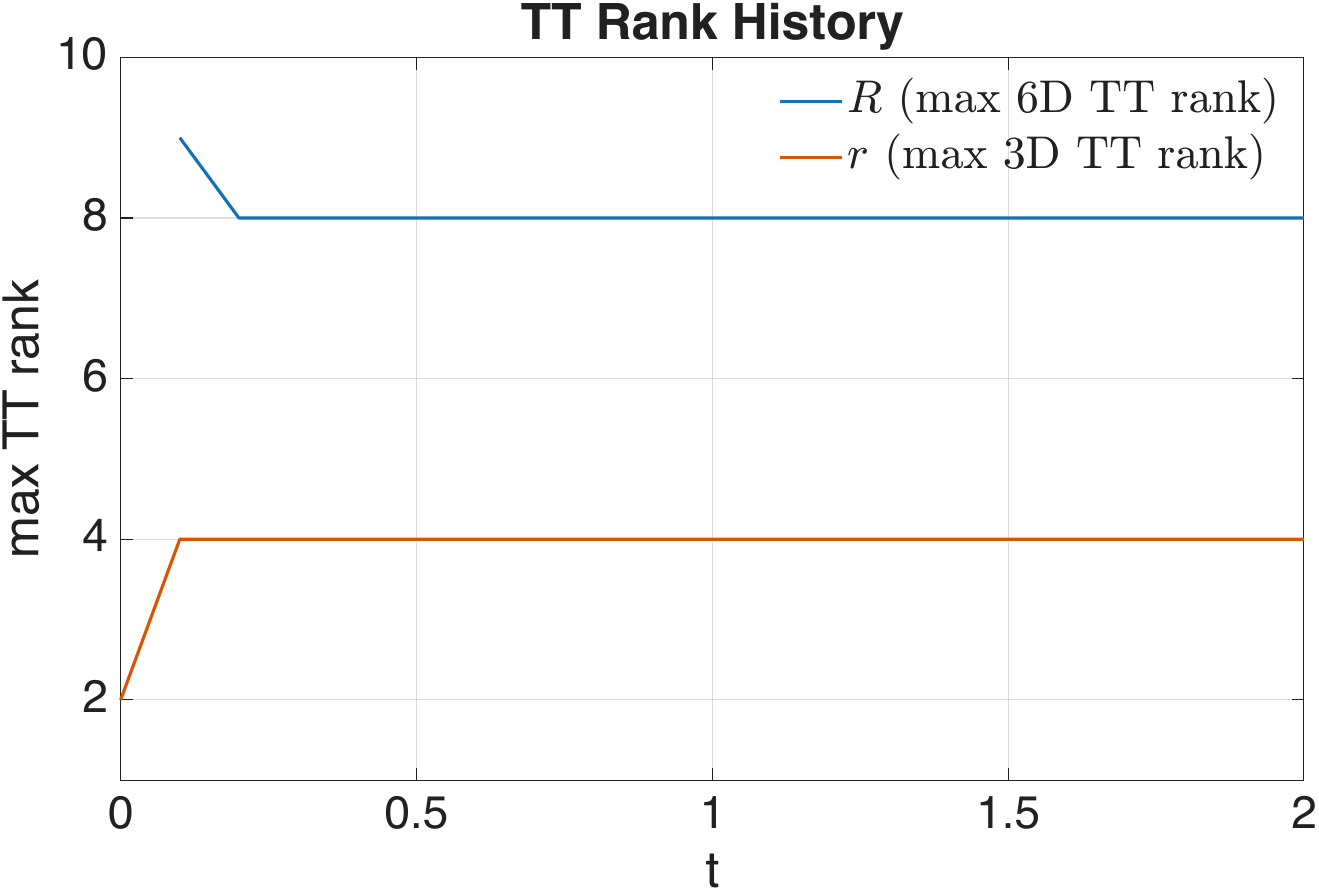}
        \caption{Maxwell molecule model with BKW initial condition. $(h,\Delta t) = (0.1, 0.1)$.}
    \end{subfigure}
    \hfill
    \begin{subfigure}[b]{0.45\textwidth}
        \centering
        \includegraphics[width=\textwidth]{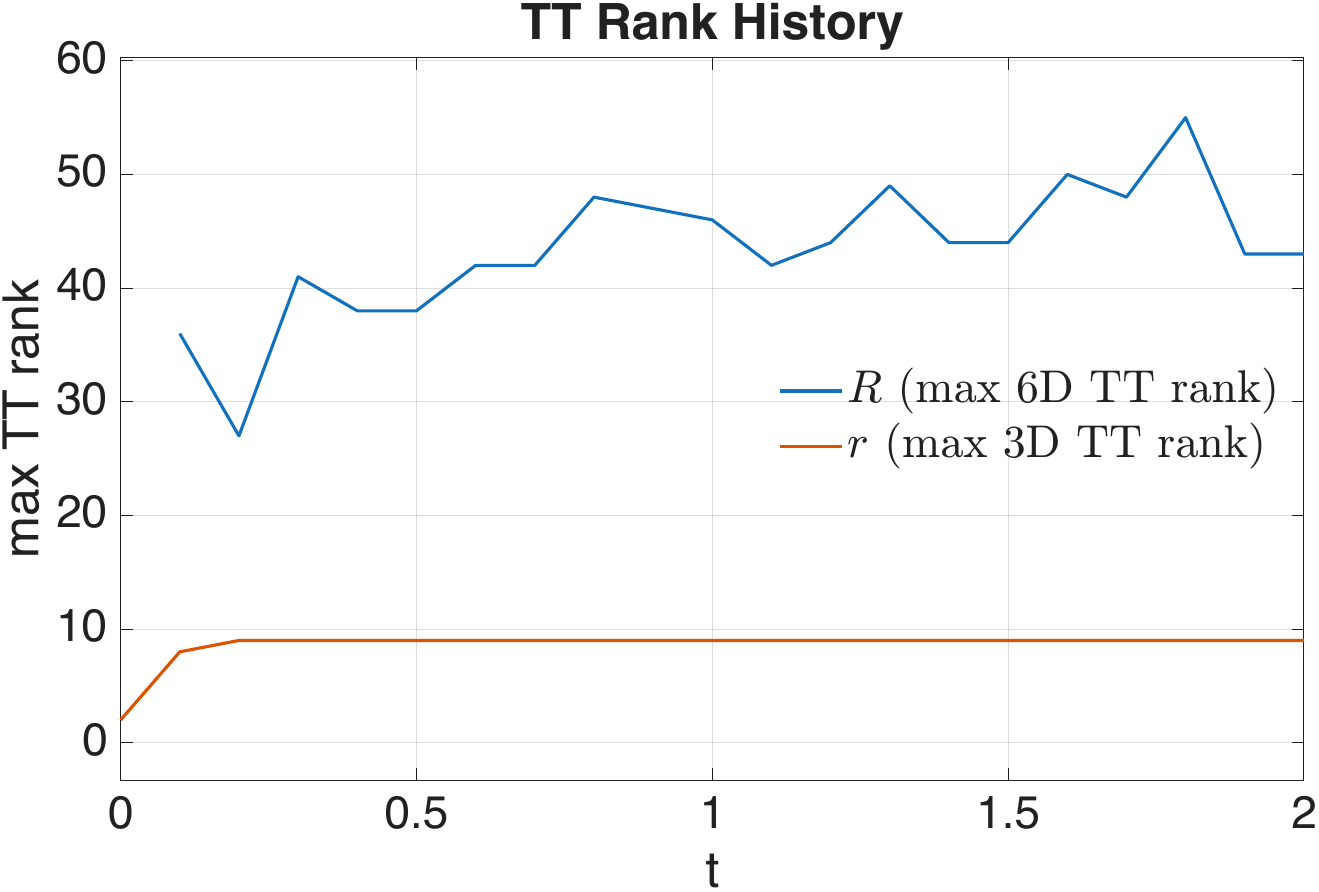}
        \caption{Hard sphere model with two Gaussians initial condition. $(h,\Delta t) = (0.1, 0.1)$. }
    \end{subfigure}
    \caption{TT rank history of the numerical solutions over time. }
    \label{fig:rank_history}
\end{figure}

Figure \ref{fig:rank_history} shows the maximal TT ranks  of the numerical solutions over time for two test cases. We can see that the TT ranks of the lifted solutions $R$ are higher than $r$ and vary with initial conditions and time.

 To show the performance of the proposed method, we let the system evolve from $t=0$ to $t=0.5$ with RK4 time discretization, where $\Delta t = 0.1$. The average computational time per time step is shown in Tables \ref{tab:time_cost_mm_cubic} and  \ref{tab:time_cost_hs_cubic} for the maxwell molecule model and hard sphere model respectively. The time cost is measured on MacBook Pro with M4 chip, with a MATLAB implementation. In particular, for tensor train operations we use the implementation provided in the GitHub repository by Oseledets et. al. \cite{tttoolbox}, with slight modifications.
 
 For fixed tolerance, the time cost is proportional to $n$, since the rank of the lifted solution is almost independent of $n$. For fixed $n$, the time cost increases with respect to the accuracy requirement, since higher accuracy requires higher rank for the lifted solution. The BKW solution with Maxwell molecule is a special case where the lifted solution is strictly low-rank, which explains the low computational cost in Table \ref{tab:time_cost_mm_cubic}. The computational results verify the complexity scaling. In particular, we note the importance of choosing appropriate tolerance according to local truncation error. If the tolerance is set to be too small, the cross approximation routine will fail to meet the tolerance and reach max number of sweep, which explains high computational cost.

\begin{table}[H]
    \centering
    \begin{tabular}{c|c|c|c}
         & $\text{tol} = 0.1\times(0.4)^{4}$ & $\text{tol} = 0.1\times(0.2)^{4}$ & $\text{tol} = 0.1\times(0.1)^{4}$\\
        $n = 32$ & 0.5 sec & * 4.5 sec & * 8.0 sec\\
        $n = 64$ & 0.8 sec & 0.9 sec  & * 33.0 sec\\
        $n = 128$ & 1.6 sec & 1.6 sec & 1.8 sec\\
    \end{tabular}
    \caption{Computational time per time step for maxwell molecule model with cubic interpolation. The "*" means that at least one cross approximation failed to meet tolerance and exited with max number of sweep.}\label{tab:time_cost_mm_cubic}
\end{table}

\begin{table}[H]
    \centering
    \begin{tabular}{c|c|c|c}
         & $\text{tol} = 0.1\times(0.4)^{4}$ & $\text{tol} = 0.1\times(0.2)^{4}$ & $\text{tol} = 0.1\times(0.1)^{4}$ \\
        $n = 32$ &  0.8 sec & 2.8 sec & * 8.2 sec \\
        $n = 64$ &  1.4 sec &  4.1 sec  & 14.5 sec \\
        $n = 128$ &  2.3 sec &  8.0 sec &  23.7 sec \\
    \end{tabular}
    \caption{Computational time per time step for hard sphere model with cubic interpolation. The "*" means that at least one cross approximation failed to meet tolerance and exited with max number of sweep.}\label{tab:time_cost_hs_cubic}
\end{table}

\subsubsection{Complexity of spectral interpolation}
Recall that with spectral interpolation, the complexity of the proposed method is $\mathcal{O}(n^{2}MR^2 r^2 X + nR^3X)$, where $n$ is the number of grid points in each direction, $M$ is the number of quadrature points on the sphere, $R$ is the max TT rank of the 6d lifted solution $F$, $r$ is the TT-rank of the 3D projected solution $f$, and $X$ is the number of sweeps. 

Following the same procedure and using the same parameter as in Section \ref{subsubsec:cubic_complexity}, the average time cost is shown in Tables \ref{tab:time_cost_mm_spectral} and \ref{tab:time_cost_hs_spectral}.  For the Maxwell molecule BKW test, the measured times are nearly independent of the tolerance because the lifted solution remains very low rank. For the hard sphere test, tighter tolerances increase the ranks and sweep counts, so the time grows with both $n$ and the requested accuracy, consistent with the rank-dependent complexity estimate. When comparing Table \ref{tab:time_cost_hs_spectral} with Table \ref{tab:time_cost_hs_cubic}, it can be observed that  for $n = 32$, the spectral solver is slightly faster than the cubic one, while in the $n=128$ case, it is much slower. This observation is also consistent with the complexity analysis: when $n$ is small, the auxiliary operation cost $\mathcal{O}(nR^3X)$ in TT cross approximation dominates, and as $n$ increases, the interpolation cost $\mathcal{O}(n^{2}MR^2 r^2 X)$ becomes the dominant term since it scales quadratically with $n$.

\begin{table}[H]
    \centering
    \begin{tabular}{c|c|c|c}
         & $\text{tol} = 0.1\times(0.4)^{4}$ & $\text{tol} = 0.1\times(0.2)^{4}$ & $\text{tol} = 0.1\times(0.1)^{4}$\\
        $n = 32$ & 0.5 sec & 0.5 sec & 0.5 sec\\
        $n = 64$ & 0.9 sec & 0.9 sec  & 0.9 sec\\
        $n = 128$ & 2.5 sec & 2.4 sec & 2.5 sec\\
    \end{tabular}
    \caption{Computational time per time step for maxwell molecule model with spectral interpolation.}
    \label{tab:time_cost_mm_spectral}
\end{table}

\begin{table}[H]
    \centering
    \begin{tabular}{c|c|c|c}
         & $\text{tol} = 0.1\times(0.4)^{4}$ & $\text{tol} = 0.1\times(0.2)^{4}$ & $\text{tol} = 0.1\times(0.1)^{4}$ \\
        $n = 32$ &  0.7 sec & 2.1 sec & * 5.7 sec \\
        $n = 64$ &  1.6 sec &  5.8 sec  & 17 sec \\
        $n = 128$ &  4.5 sec &  14.1 sec &  44.1 sec \\
    \end{tabular}
    \caption{Computational time per time step for hard sphere model with spectral interpolation. The ``*" means that at least one cross approximation had exited above tolerance.}
    \label{tab:time_cost_hs_spectral}
\end{table}

\subsection{Necessity of conservation correction}
 In prior subsections, all the low-rank numerical results are obtained with conservation correction.  In this subsection, we run the numerical tests with and without conservation correction to validate its effect. Using cubic interpolation and RK4, we set $L=6.4$, $M=74$ Lebedev points, $n=64$, $\Delta t=0.2$, a maximum number of greedy sweeps equal to the number of grid points per velocity direction, and TT-rounding tolerance equal to one tenth of the cross-approximation tolerance.

\begin{figure}[H]
    \centering
    \begin{subfigure}[b]{0.45\textwidth}
        \centering
        \includegraphics[width=\textwidth]{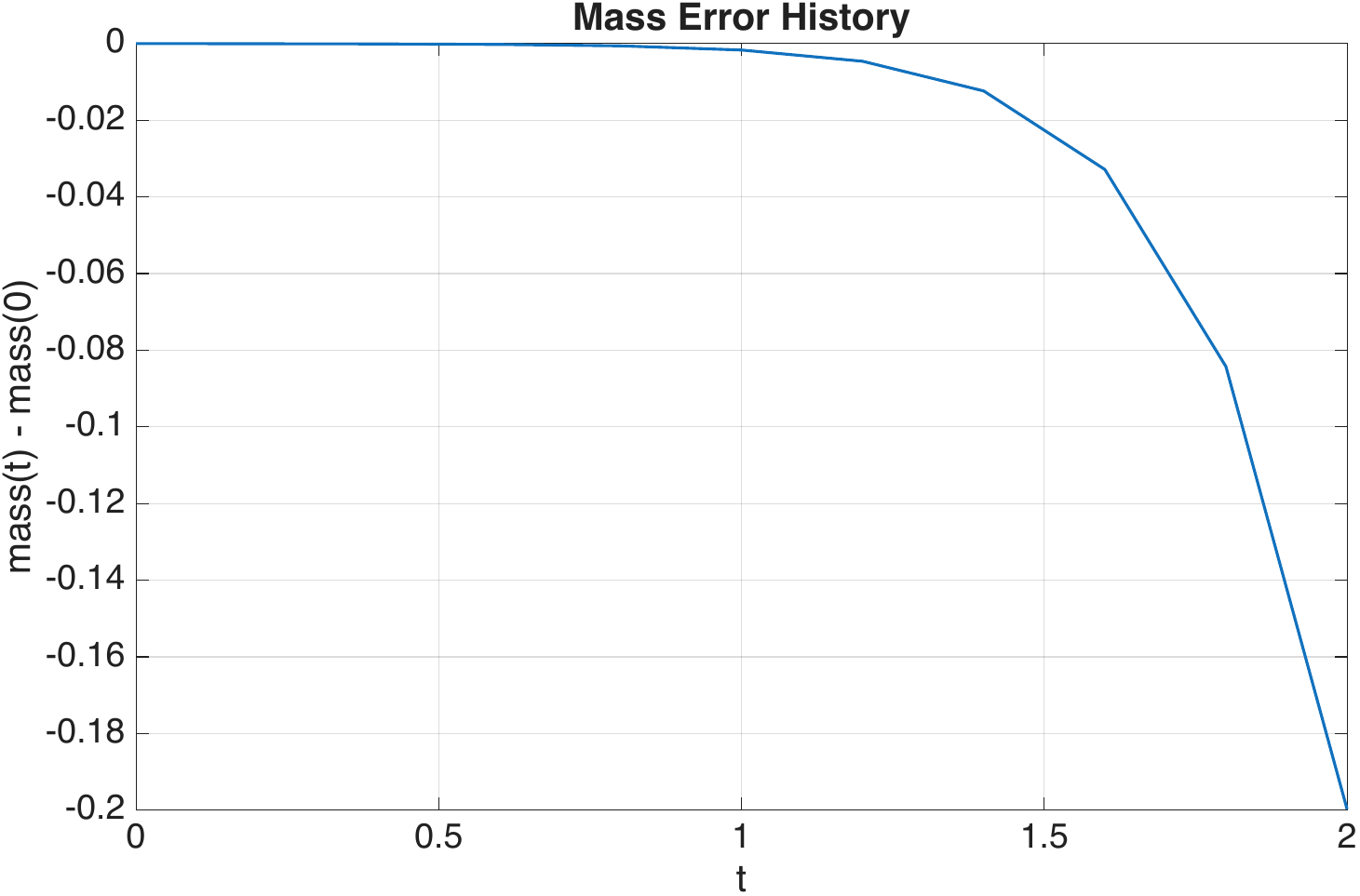}
    \end{subfigure}
    \hfill
    \begin{subfigure}[b]{0.45\textwidth}
        \centering
        \includegraphics[width=\textwidth]{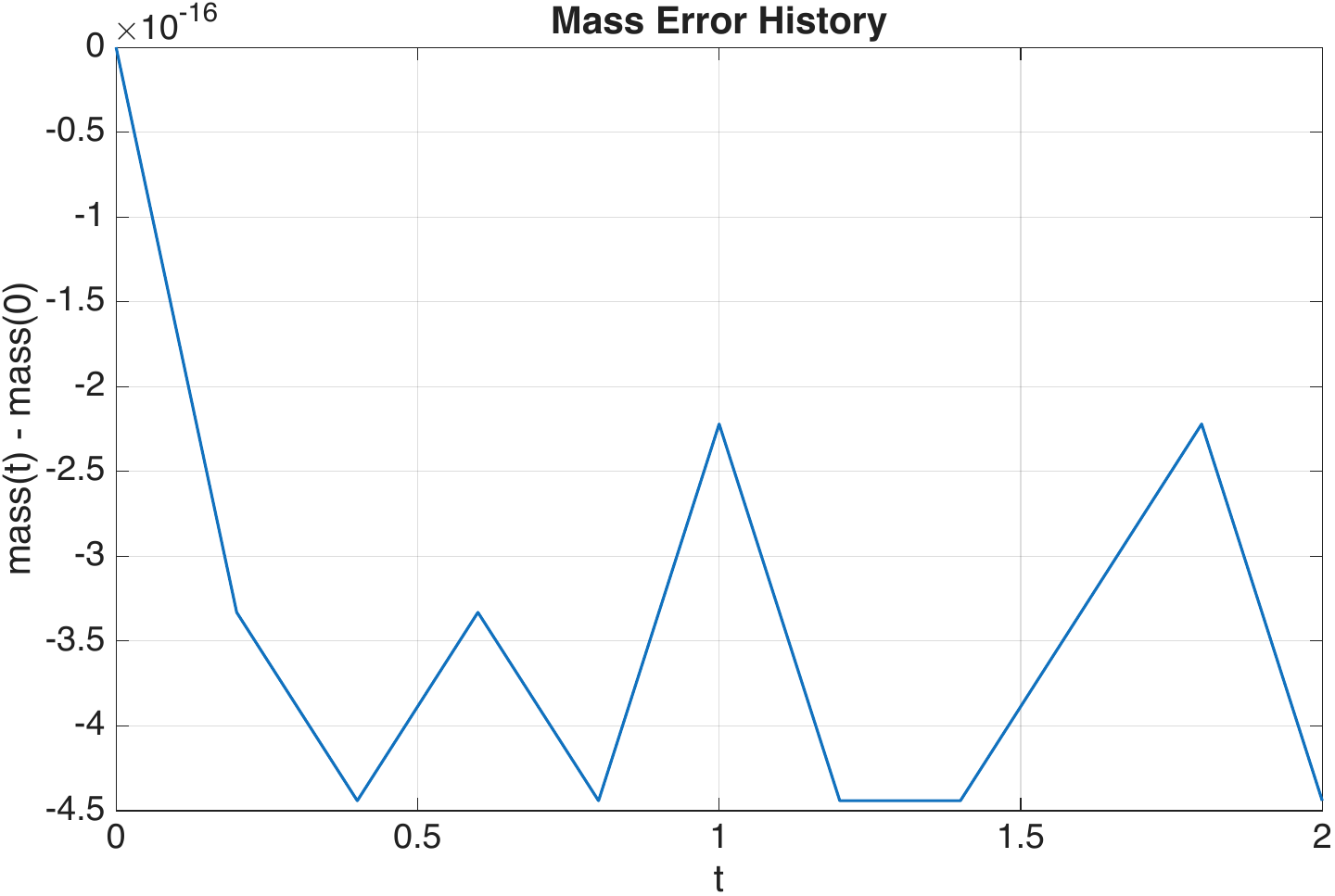}
    \end{subfigure}
    \begin{subfigure}[b]{0.45\textwidth}
        \centering
        \includegraphics[width=\textwidth]{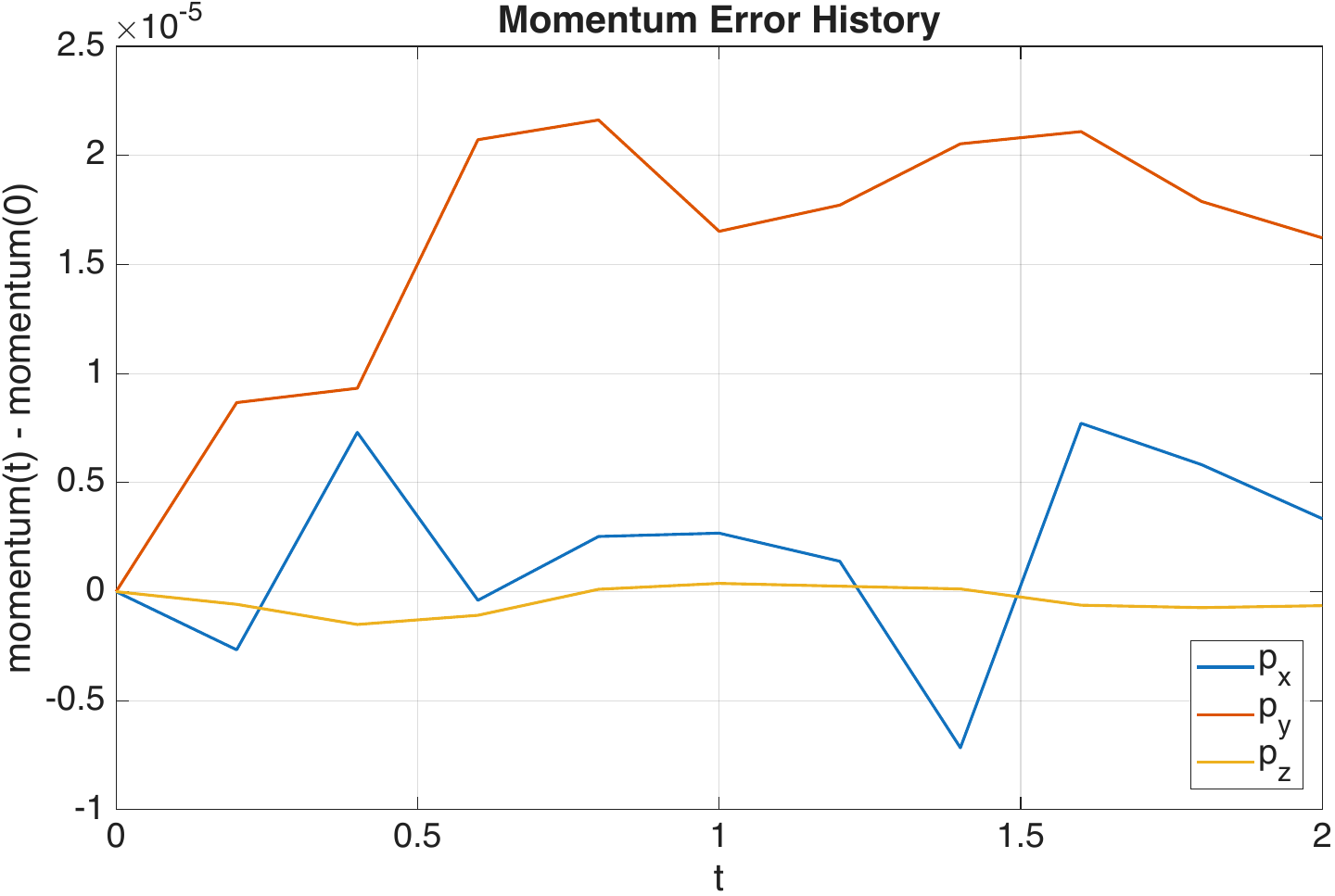}
    \end{subfigure}
    \hfill
    \begin{subfigure}[b]{0.45\textwidth}
        \centering
        \includegraphics[width=\textwidth]{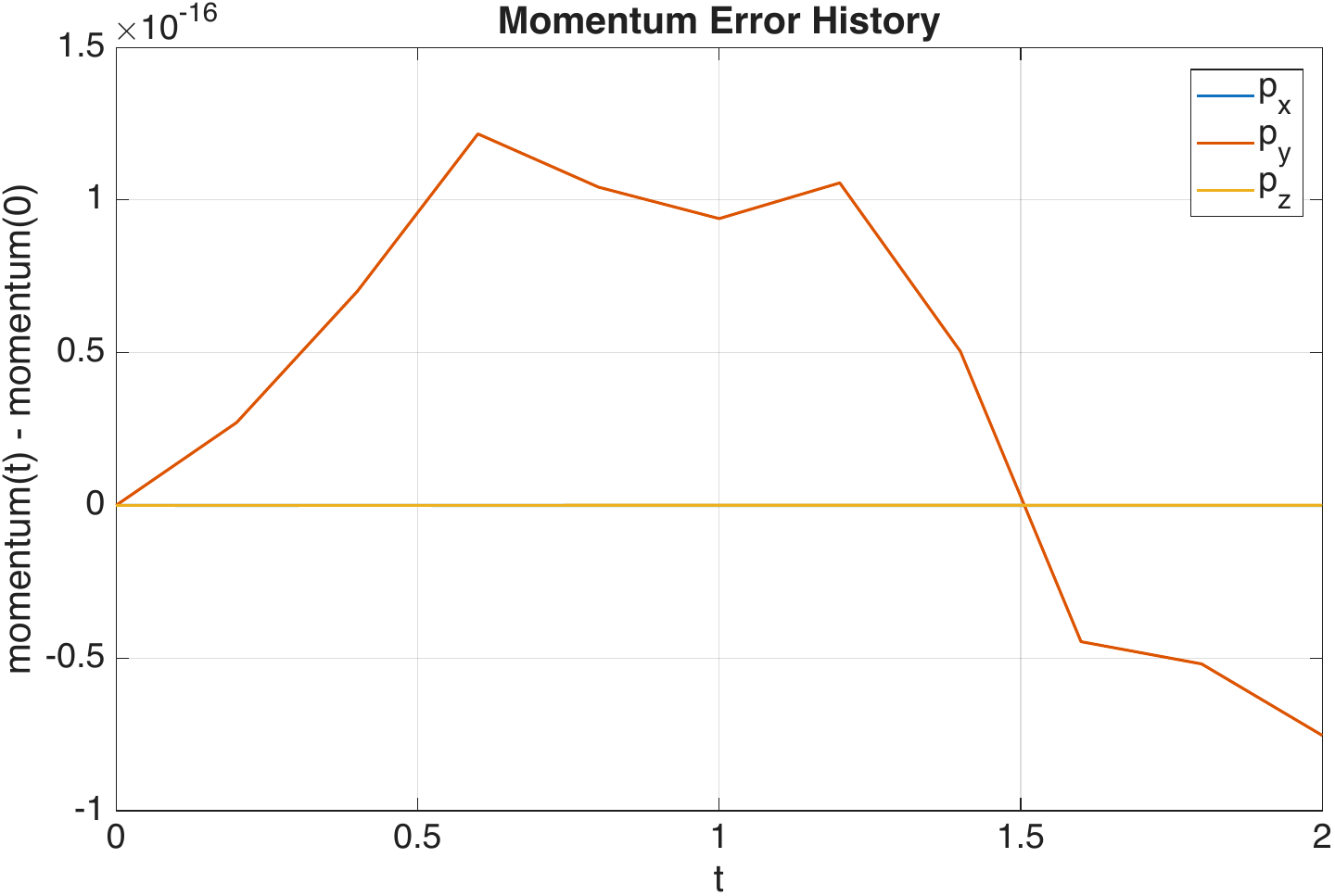}
    \end{subfigure}
    \begin{subfigure}[b]{0.45\textwidth}
        \centering
        \includegraphics[width=\textwidth]{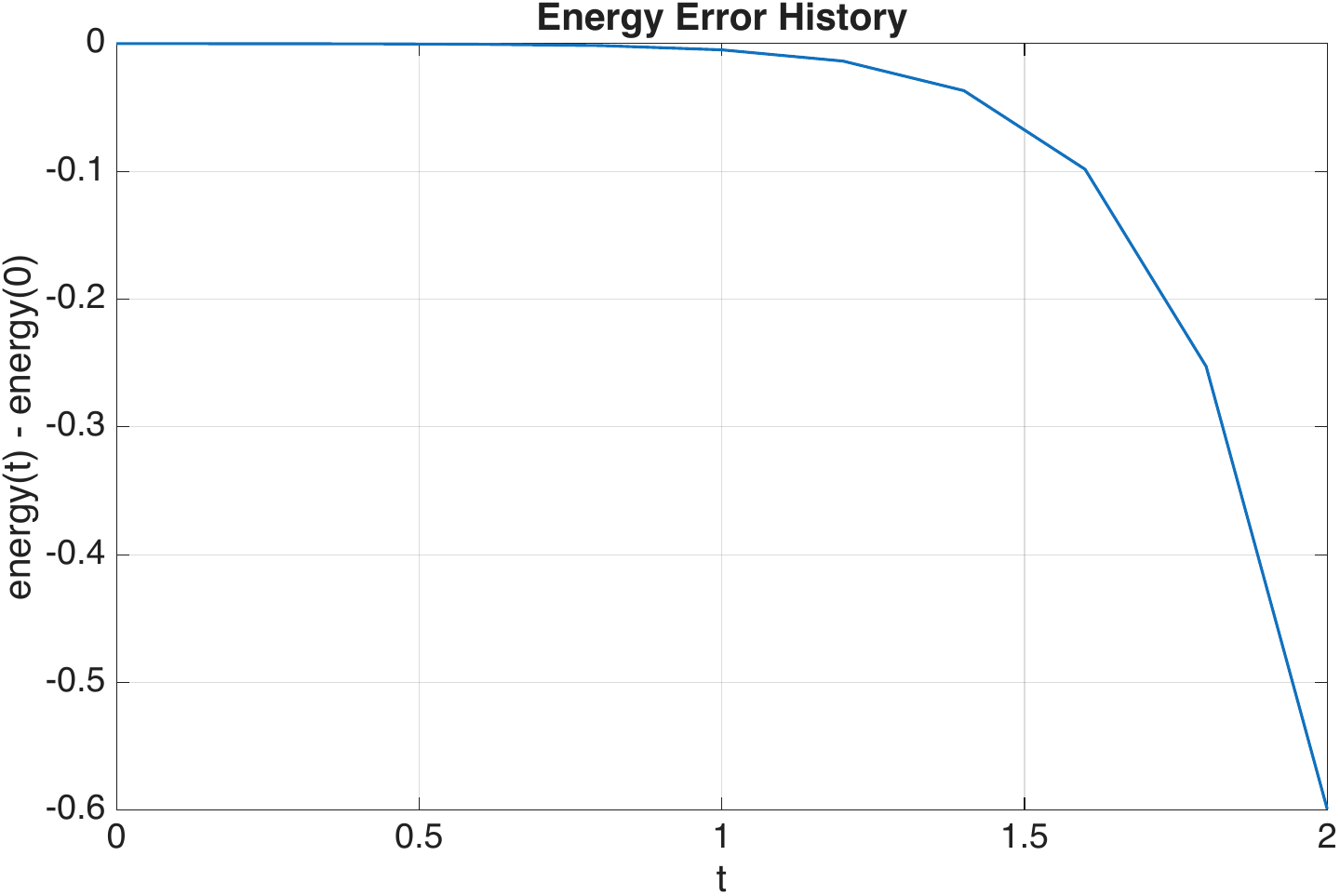}
    \end{subfigure}
    \hfill
    \begin{subfigure}[b]{0.45\textwidth}
        \centering
        \includegraphics[width=\textwidth]{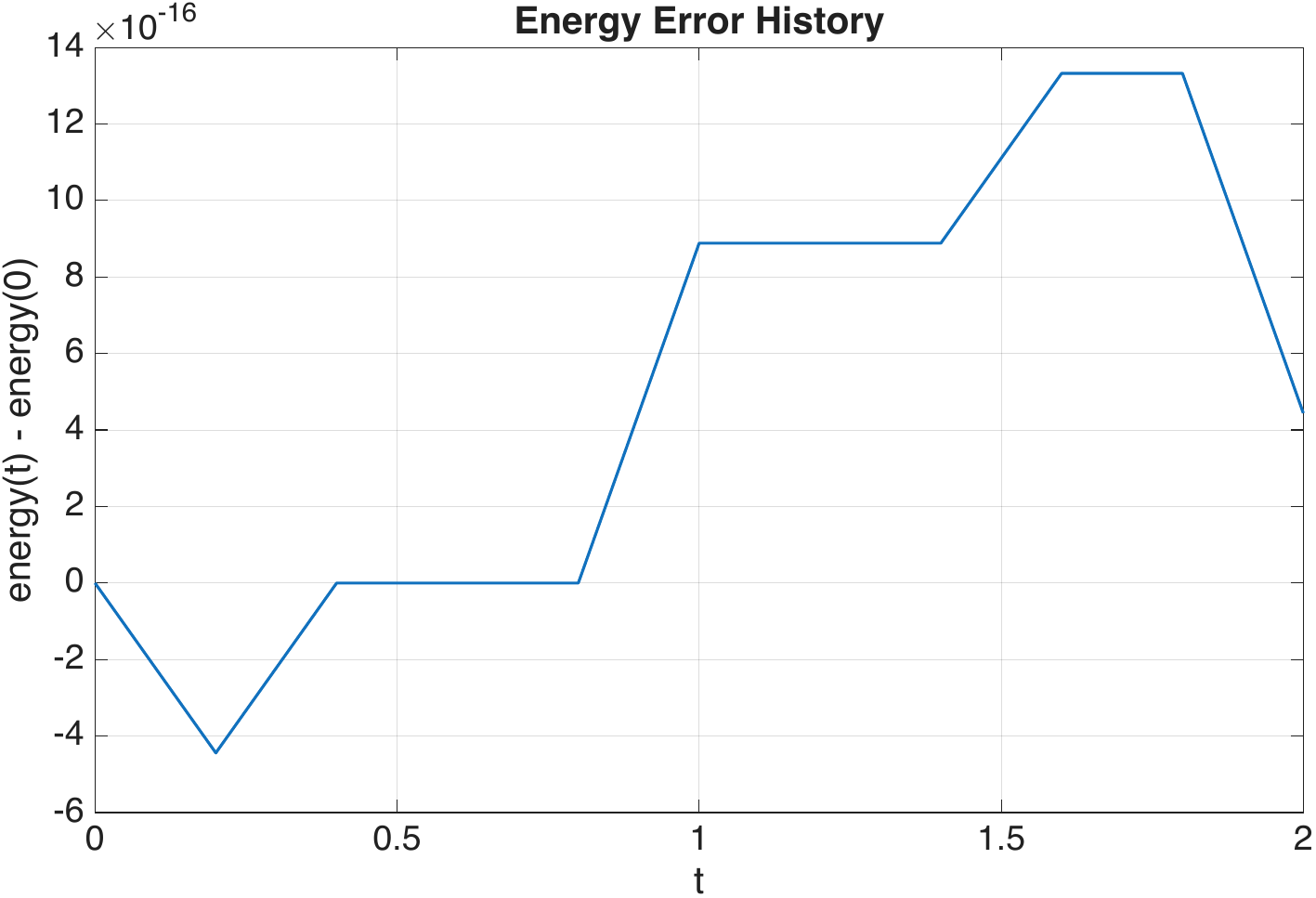}
    \end{subfigure}
    \caption{Numerical values of mass, momentum and energy. Left column: without conservation correction; right column: with conservation correction.}
    \label{fig:conservation}
\end{figure}

\begin{figure}[H]
    \centering
    \begin{subfigure}[b]{0.45\textwidth}
        \centering
        \includegraphics[width=\textwidth]{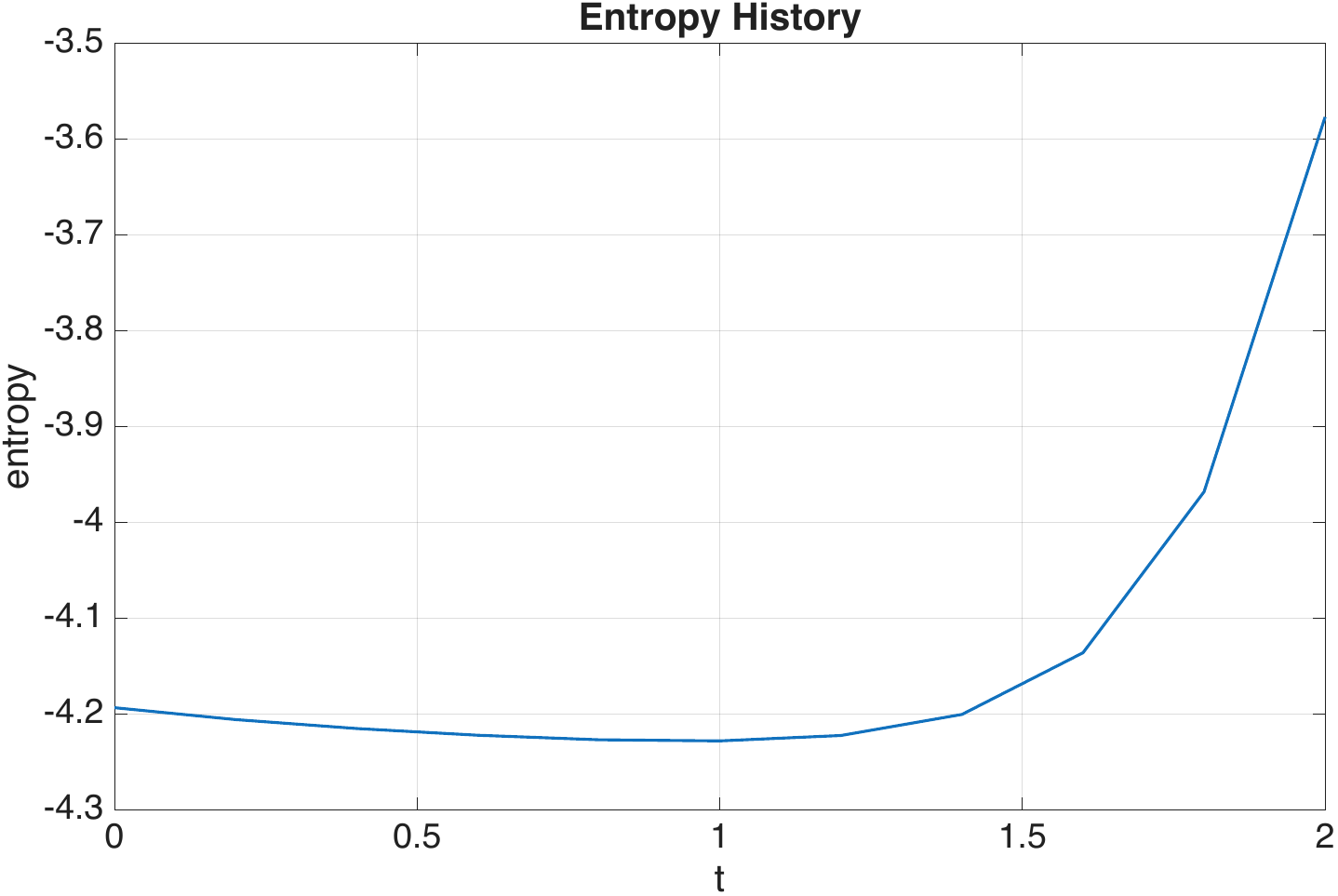}
    \end{subfigure}
    \hfill
    \begin{subfigure}[b]{0.45\textwidth}
        \centering
        \includegraphics[width=\textwidth]{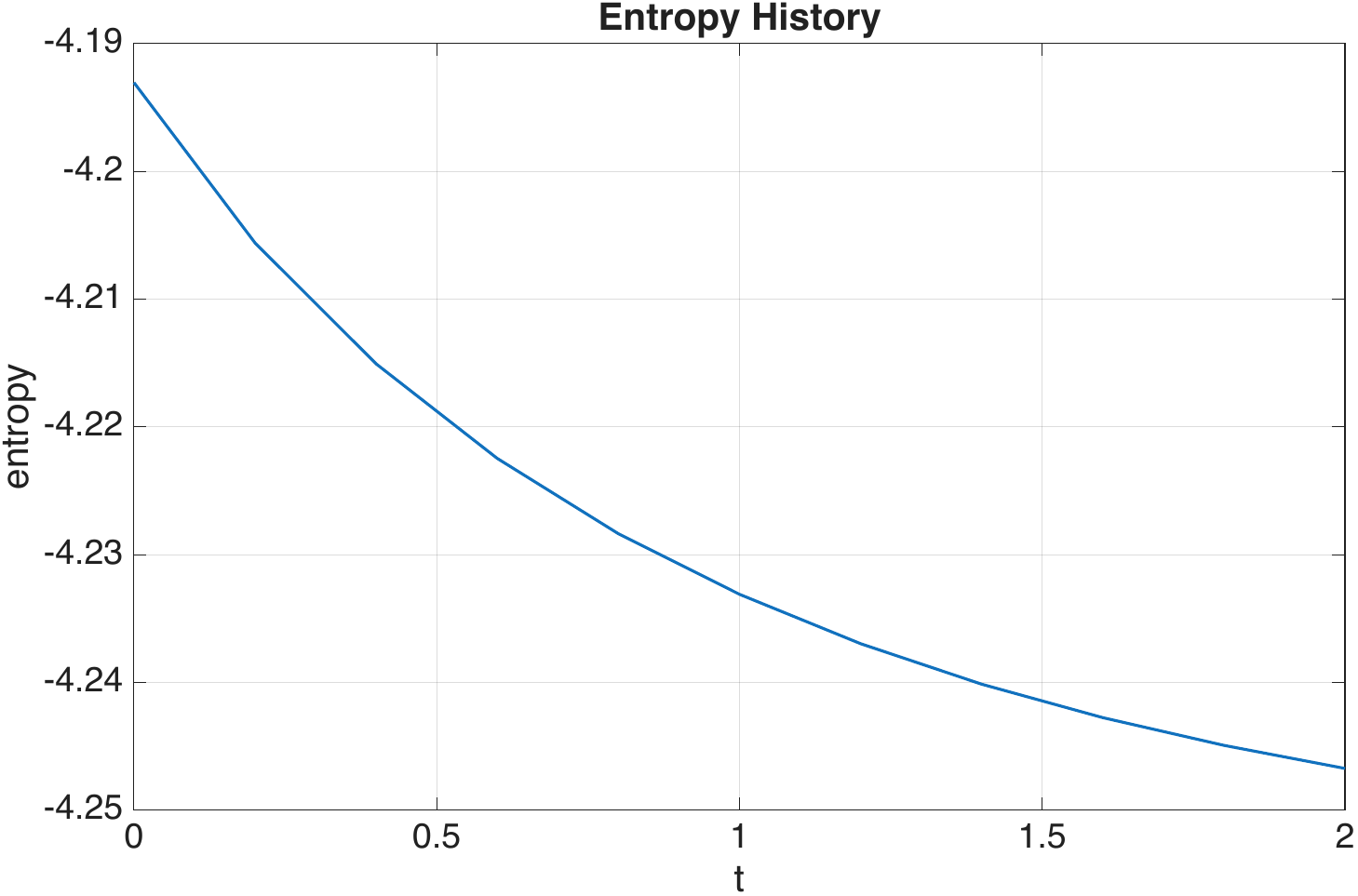}
    \end{subfigure}
    \caption{Numerical entropy. Left column: without conservation correction; right column: with conservation correction.}
    \label{fig:entropy}
\end{figure}

The numerical results are gathered in  Figures \ref{fig:conservation} and \ref{fig:entropy}. It is clear that the numerical results without conservation correction present large errors in the moments. Indeed, an essential condition for the lifting-projection procedure to work is that the 3D solution $f$ should be a probability distribution s.t. $\int_{\mathbb{R}^{3}} f(\mathbf{v}) d \mathbf{v} = 1$. More specifically, in Algorithm \ref{alg:rk4}, when computing 
\begin{equation*}
    q(f^{(s)})= \Pi Q (f_{h}^{(s)}\otimes f_{h}^{(s)}) = \frac{\Pi (I + Q \Delta t) (f_{h}^{(s)}\otimes f_{h}^{(s)}) - f_{h}^{(s)}}{\Delta t}, 
\end{equation*}
we have implicitly assumed that
\begin{equation*}
    \Pi (f_{h}^{(s)}\otimes f_{h}^{(s)}) - f_{h}^{(s)} =0,
\end{equation*}
which is only true when mass conservation is preserved in every step. The computational results validate the essential roles of conservation corrections.

\section{Conclusion} \label{sec:conclusion}
We  presented a novel conservative, low-rank TT Boltzmann solver consisting of two independent
  contributions. The first is the lifting-projection (LP) scheme, which is a computational framework motivated by recent theoretical breakthrough \cite{guillen2025landau,imbert2026monotonicity,guillen2025landau2}  that lifts the
  nonlinear 3D Boltzmann equation to a 6D linear Kac master equation, evolve for one time step, and projects the result back
  to 3D. 
  The second contribution is a fast low-rank tensor method. We represent the lifted solution in TT format and compute 
  it via a TT cross algorithm.  With cubic interpolation, the cost of the proposed method is $\mathcal{O}(nMR^{2}r^{2}X + n R^3 X)$, where $n$ is the number of grid points in each direction, $M$ is the number of quadrature points on the sphere, $R$ is the max TT rank of the 6d lifted solution $F$, $r$ is the TT-rank of the 3D projected solution $f$, and $X$ is the number of sweeps. With spectral interpolation, the cost is $\mathcal{O}(n^{2}MR^{2}r^{2}X)$. When the solution exhibits low-rank property, this method, which scales linearly with $n,$ offers huge computational savings compared to existing methods in the literature. We further propose a simple conservation correction procedure that is compatible with TT format that ensures mass, momentum and energy conservation. 

  We perform extensive numerical tests, benchmarking the performance of the methods in terms of accuracy and computational efficiency. The numerical experiments demonstrate that the $xyz$-ordering in TT leads to significantly lower  ranks for the
  lifted solution compared to the $vw$-ordering. Convergence tests show that the methods achieve the expected accuracy and numerical conservations are maintained.

  Several directions remain open for future work. The method currently handles
  only the spatially homogeneous Boltzmann equation. Extending it to the full
  spatially inhomogeneous case is a natural next step.  We also plan to investigate applications to Landau equation and study of the asymptotic preserving schemes.

\section*{Acknowledgements}
This material is based upon work supported by the National Science Foundation under Grant No.
DMS-2424139 while the  authors were in residence at the Simons Laufer Mathematical Sciences
Institute in Berkeley, California, during the Fall 2025 semester.

\bibliographystyle{plain}
\bibliography{main.bib}

\end{document}